\documentclass[12pt]{amsart}

\usepackage{amssymb,euscript,array,amscd}

\newtheorem{Lemma}{Lemma}[section]\newcommand{\bel}{\begin{Lemma}}\newcommand{\eel}{\end{Lemma}}
\newtheorem{Proposition}[Lemma]{Proposition}\newcommand{\bprop}{\begin{Proposition}}\newcommand{\eprop}{\end{Proposition}}
\newtheorem{Theorem}[Lemma]{Theorem}\newcommand{\bthe}{\begin{Theorem}}\newcommand{\ethe}{\end{Theorem}}

\newcommand{\bpr}{~\\{\em Proof.~}}\def\epr{$\bull$\\}
\newtheorem{Remark}[Lemma]{Remark}\newcommand{\beR}{\begin{Remark}\rm}\newcommand{\eeR}{\end{Remark}}
\newtheorem{Definition}[Lemma]{Definition}\newcommand{\bed}{\begin{Definition}}\newcommand{\eed}{\end{Definition}}
\newtheorem{Example}[Lemma]{Example}\newcommand{\bex}{\begin{Example}\rm}\newcommand{\eex}{\end{Example}}
\newtheorem{Corollary}[Lemma]{Corollary}\newcommand{\bcor}{\begin{Corollary}}\newcommand{\ecor}{\end{Corollary}}
\newtheorem{DefProp}[Lemma]{Definition-Proposition}\newcommand{\bdp}{\begin{DefProp}}\newcommand{\edp}{\end{DefProp}}
\newtheorem{Property}[Lemma]{Property}\newcommand{\bpro}{\begin{Property}}\newcommand{\epro}{\end{Property}}

\newcommand{\beq}{\begin{equation}}\newcommand{\eeq}{\end{equation}}
\newcommand{\bem}{\begin{displaymath}}\newcommand{\eem}{\end{displaymath}}
\newcommand{\beqa}{\begin{eqnarray}}\newcommand{\eeqa}{\end{eqnarray}}
\newcommand{\bee}{\begin{enumerate}}\newcommand{\eee}{\end{enumerate}}
\newcommand{\bei}{\begin{itemize}}\newcommand{\eei}{\end{itemize}}
\newcommand{\bet}{\begin{tabular}{cccccccc}}\newcommand{\eet}{\end{tabular}}
\newcommand{\bpm}{\begin{pmatrix}}\newcommand{\epm}{\end{pmatrix}}
\newcommand{\bM}{\begin{matrix}}\newcommand{\eM}{\end{matrix}}
\newcommand{\ber}{\begin{array}{l}}\newcommand{\eer}{\end{array}}

\def\bull{\vrule height .9ex width .9ex depth -.1ex }

\def\bl{\langle}\def\br{\rangle}
\def\da{\downarrow}\def\isom{\xrightarrow{\sim}}
\def\into{\stackrel{i}{\hookrightarrow}}\def\norm{\stackrel{\nu}{\to}}

\newcommand{\quotient}[2]{{\left.\raisebox{0.4ex}{$#1$}\!\!\middle/\!\!\raisebox{-0.4ex}{$#2$}\right.}}

\newcommand{\li}{~\\ $\bullet$ }

\def\cM{\mathcal{M}}\def\cN{\mathcal{N}}\def\cO{\mathcal{O}}

\def\cm{{\frak m}}

\def\one{{1\hspace{-0.1cm}\rm I}}\newcommand{\zero}{0\put(-2.5,0.1){\line(0,1){7.2}}}

\def\k{k}
\def\P{\mathbb{P}}
\def\R{\mathbb{R}}

\def\al{\alpha}\def\be{\beta}

\def\tA{\tilde{A}}\def\tC{{\tilde{C}}}

  \def\tX{{\tilde{X}}}

\def\prodl{\prod\limits}

\def\cMv{{\cM^{\vee}}}\def\cNv{{\cN^{\vee}}}

\def\smin{\setminus}\def\sset{\subset}

\def\dr{determinantal representation}\def\Dr{Determinantal representation}
\def\mf{matrix factorization}
\def\mg{maximally generated}
\def\XXS{$X'/X$-saturated }
\def\CM{Cohen-Macaulay }\def\cOn{\cO_{(k^n,0)}}\def\cOx{\cO_{(X,0)}}

\def\omp{ordinary multiple point}

\title[D\MakeLowercase{eterminantal representations}]{D\MakeLowercase{ecomposability of local \dr s of hypersurfaces}}
\author[D.K\MakeLowercase{erner}]{D\MakeLowercase{mitry} K\MakeLowercase{erner}}
\address{Department of Mathematics, University of Toronto, 40 St. George Street, Toronto, Canada.}
\email{dmitry.kerner@utoronto.ca}
\author[V.V\MakeLowercase{innikov}]{V\MakeLowercase{ictor} V\MakeLowercase{innikov}}
\address{Department of Mathematics, Ben Gurion University of the Negev, P.O.B. 653, Be'er Sheva 84105, Israel.}
\email{vinnikov@math.bgu.ac.il}
\thanks{
Part of the work was done during postdoctoral stay of D.K. in Mathematics Department of Ben Gurion
University, Israel.  Both authors were supported by the Israel Science Foundation.
\\ The authors thank G.Belitski, R.O.Buchweitz, I.Burban, G.M.Greuel, M.Leyenson and I.Tyomkin
for numerous important discussions.
}
\date{\today}

\begin{document}
\maketitle  \setcounter{secnumdepth}{4}\setcounter{tocdepth}{1}
\begin{abstract}
Let $\cM$ be a matrix whose entries are power series in several variables and determinant $\det(\cM)$ does not
vanish identically. The equation $\det(\cM)=0$ defines a hypersurface singularity and the (co)-kernel
of $\cM$ is a maximally Cohen-Macaulay module over the local ring of this singularity.

Suppose the determinant $\det(\cM)$ is reducible, i.e. the hypersurface is locally reducible.
A natural question is whether the matrix is equivalent to a block-diagonal or at least to an upper-block-triangular.
 (Or whether the corresponding module is decomposable or at least is an extension.)
We give various necessary and sufficient criteria.

Two classes of such matrices of functions appear naturally in the study of decomposability:
those with many generators (e.g. \mg\ or Ulrich maximal) and those that descend from
birational modifications of the hypersurface by pushforwards (i.e. correspond to modules over
bigger rings).  Their properties are studied.
\end{abstract}
\tableofcontents

\section{Introduction}\label{Sec.Introduction}
Let $\k$ be an algebraically closed, normed, complete field of characteristic zero, e.g. the complex numbers.
Let $\k^n$ be the affine space of dimension $n$ over $k$, let $(\k^n,0)$ be a small neighborhood of the origin.
Denote the corresponding ring of regular functions by $\cOn$. This means:
\\\ - rational functions that are regular at the origin: $k[x_1,..,x_n]_{(\cm)}$,
\\\ - or locally converging power series, $k\{x_1,..,x_n\}$,
\\\ - or formal power series, $k[[x_1,..,x_n]]$.

In this paper by a curve/hypersurface we always mean the germ at the singular point,
which is assumed to be the origin $0\in(\k^n,0)$. The hypersurfaces are considered with their multiplicities,
not just as zero sets of functions. We denote the zero/identity matrices by $\zero$, respectively $\one$.
Denote by $\cMv$ the adjoint matrix of $\cM$, so $\cM\cMv=\det(\cM)\one$.

\subsection{Setup}
Let $\cM$ be a $d\times d$ matrix whose entries are functions in $\cOn$.
We always assume \mbox{$f=\det(\cM)\not\equiv0$} and $d>1$, and the matrix vanishes at the origin, $\cM|_0=\zero$.
So the matrix defines the (algebraic, analytic or formal) {\em hypersurface}
 $(X,0)=\{\det(\cM)=0\}\sset(\k^n,0)$, of dimension $(n-1)$. This hypersurface is mostly singular, can be reducible/non-reduced.

Such matrices of functions occur in various fields. In algebraic geometry they are called \underline{local
\dr s} of the hypersurface $(X,0)$. They correspond also to some elements of the local
class group $Cl(X,0)$. In commutative algebra they are appear as matrix factorizations or the syzygy matrices in the
resolutions of modules over  hypersurface singularities. They appear as
local homomorphisms of vector bundles (and degeneracy loci), as maps $(\k^n,0)\to Mat(d\times d)$,
as matrix families in operator theory, as transfer functions in control theory,
in semi-definite programming etc.
For a short mixture of results cf. \S\ref{Sec.Intro Brief Overview}.
\\\

With the applications in mind, the \dr s are studied up to the {\em local equivalence}
 $\cM\sim A\cM B$ for $A,B\in GL(d,\cOn)$, i.e. $A,B$ are invertible at the origin.
 This equivalence preserves the embedded hypersurface pointwise. Any matrix is locally equivalent to
 a block-diagonal, $\cM\sim\one\oplus\cM'$, where $\cM'|_0=\zero$, property \ref{Thm.Localiz.Chip.off.Unity}.
 Hence we mostly assume that $\cM$ vanishes at the origin.

\Dr s are well studied in simple cases, e.g. when the hypersurface singularity is of one of $A,D,E$ types, or
for locally irreducible plane curve singularities. For more complicated singularities the \dr s are not
so well understood. In particular the following decomposability question has not been addressed.

{\bf Problem} {\em Suppose $det(\cM)$ is reducible, i.e. the hypersurface $\{\det\cM=0\}$ is locally decomposable. When
is $\cM$ decomposable, i.e. $\cM\sim\cM_1\oplus\cM_2$? When is $\cM$ an extension, i.e. equivalent to
an upper-block-triangular matrix?}

In more modern language this question reads: when is the (co)kernel module
 decomposable or extension?

The main goal of this work is to treat this question and provide necessary and/or sufficient
decomposability criteria. Our analysis shows two particular types of matrix functions (or modules)
with nice properties: those of "maximal size"
and those that "descend from a modification of singularity". They admit especially strong criteria.
On the other hand, they happen to be particularly
important in applications, \cite{Livšic-Kravitsky-Markus-Vinnikov-book}. Any such decomposability
criterion is useful as it reduces a "matrix problem" to a simpler blocks.
\subsection{Maximally generated \dr s, i.e. matrices of maximal size}
Restrict $\cM$ to the hypersurface $(X,0)$, i.e. consider $\cM$ as a matrix with entries in the quotient
ring $\cOx=\cOn/(\det\cM)$.
Two \dr s are equivalent over $\cOn$ iff they
are equivalent over $\cOx$, (proposition \ref{Thm.Equivalence.Over.Hypersurface.Implies.Ordinary}).

Let $mult(X,pt)$ be the multiplicity of the hypersurface at the point, i.t. the order of vanishing of $\det\cM$.
At each point $corank\cM|_{pt}\le mult(X,pt)$, (property \ref{Thm.Localiz.Chip.off.Unity}).
This motivates the following
\bed\label{Def Maximality and Weak Maximality}
The representation is called \underline{\mg}  at the point $0\in (X,0)\sset(\k^n,0)$
 if $corank\cM|_{0}=mult(X,0)$.  The representation is called \mg\ \underline{near} the point $0\in X\sset\k^n$
 if it is \mg\ in some neighborhood of $0\in\k^n$.
\eed
For example, any \dr\ of $X$ is \mg\ at any smooth point of $X$.
\Dr s that are \mg\ at the origin $0\in\k^n$ correspond to Ulrich-maximal Cohen-Macaulay modules, \cite{Ulrich84}.
For an isolated hypersurface singularity (e.g. reduced plane curves) being \mg\ at the point and \mg\ near
the point is the same.

Sometimes we specify the neighborhood where $\cM$ is \mg. For example it can be the set of all the
smooth points of $X$ (or the smooth points of the reduced locus $X_{red}$), so the neighborhood is punctured.
Or the set of all the points where the multiplicity of $X$ is bounded by some number.

The notions of
\mg\ at the point and \mg\ in the punctured neighborhood of the point are essentially different. For example,
if $(X,0)$ is an isolated singularity and $\cM$ any of its \dr s then the block-diagonal matrix $\oplus^r\cM$
is a \dr\ of $(rX,0)$,
\mg\ on the punctured neighborhood of the origin. On the other hand, if $\cM$ is a \mg\ \dr\ of $(X,0)$, then
by inserting into the matrix $\oplus^r\cM$, above the diagonal blocks, some generic polynomials, vanishing
at the origin, we get a presentation which is \mg\ at the origin but not in the punctured neighborhood of the origin.
\subsection{Saturated \dr s}
Let $(X',0)\stackrel{\nu}{\to}(X,0)$ be a \underline{\em finite modification}. Here $(X',0)$ is a multi-germ and $\nu$ is
 a proper, surjective, finite, birational morphism of pure dimensional schemes that is an isomorphism outside the
 singular locus of $(X,0)$.
 For example, it could be (in the trivial case) the identity isomorphism: $(X,0)\isom(X,0)$.
 Or, in the 'maximal case',
  the normalization $(\tX,0)\to(X,0)$. For the details cf.
  \S\ref{Sec.Background.Normalization.Intermediate.Modifications}.

Associated to this morphism is the relative adjoint ideal $Adj_{X'/X}\sset\cOn$, defined by
\beq
Adj_{X'/X}:=\{g\in\cOn:\ \ \nu^*(g)\cO_{(X',0)}\sset\nu^{-1}\cOx,\text{ or alternatively }g\nu_*\cO_{(X',0)}\sset\cOx\}
\eeq
Note that $X$ can be non-reduced here. The restriction of $Adj_{X'/X}$ to $(X,0)$ defines the
 relative conductor ideal $I^{cd}_{X'/X}=Ann_{\cOx}\quotient{\cO_{(X',0)}}{\cOx}\sset\cOx$. For more detail cf.
\S\ref{Sec.Background.Adjoint.Conduct.Ideal}.

\bed\label{Def.Saturated.Det.Reps}
 The \dr\ $\cM$ is called \underline{$X'/X$-saturated} if every element of the adjoint matrix
 $\cMv$ belongs to the adjoint ideal $Adj_{X'/X}$.
\eed
(This definition is easy to check in particular cases. In \ref{Sec.Intro.Contents.Saturated.Modules} we give
an equivalent but more conceptual definition.)
Note that the properties of being \mg\ or $X'/X$ saturated are invariant with respect to the local equivalence.

\subsection{Contents of the paper} The matrices of functions appear in various fields, our results are
 relevant for broad audience. Hence we describe briefly several approaches
in \S\ref{Sec.Intro Brief Overview} and throughout the paper we recall some known facts.
We provide many (counter-)examples. In \S\ref{Sec.Applications} we give some applications of
the decomposability to the study of \dr s of particular hypersurface singularities.
\\
\\
Section \ref{Sec.Background} contains preliminaries and background.
In \S\ref{Sec.Background.Base.Ring} we discuss various notions of locality, i.e. the dependence on
the base rings: $\k[]_{(\cm)}$, $\k\{\}$ or $\k[[]]$.
In \S\ref{Sec.Background.Hypersurface.Singularities} we discuss curve and hypersurface singularities
and their finite modifications $(X',0)\to(X,0)$. In \S\ref{Sec.Background.Adjoint.Conduct.Ideal} we discuss
the corresponding conductor and adjoint ideals, $I^{cd}_{X'/X}$ and $Adj_{X'/X}$.
In \S\ref{Sec.Background.Kernel.Modules} we introduce the \underline{(co)kernel of a} \underline{\dr}, $E$, which is a
maximal \CM module of rank
1 over $(X,0)$.

\subsubsection{Decomposability of modules with many generators.}
An arbitrary \dr~ cannot be brought to an upper-block-triangular form,
even in the case of plane curve singularity that is an \omp, i.e. the union of smooth pairwise
non-tangent branches. However, modules with many generators tend to be decomposable or extensions.

{\em Proposition \ref{Thm.Decomposability.Max.Gener.Can be Brought to Upper Triangular}.
Suppose $E_{(X,0)}$ is minimally generated by $d(E)$ elements and its restrictions $E_i:=E|_{(X_i,0)}/Torsion$
are generated by $d(E_i)$ elements. Suppose the restriction $tr(E)_i$, (of Auslander's transpose) is
generated by $d(tr(E)_i)$ elements.
Then $E$ is an extension iff $d(E)=d(E_1)+d(tr(E)_2)$ or $d(E)=d(E_2)+d(tr(E)_1)$.}

In other words, suppose $\cM_{d\times d}$ vanishes at the origin. Consider the $\cO_{(X_i,0)}$ module
spanned by the columns of $\cMv|_{(X_i,0)}$. Suppose this module is minimally generated by $d(E_i)$ columns.
Suppose the module of the rows of $\cMv|_{(X_i,0)}$ is generated by $d(tr(E)_i)$ rows.
 Then $\cM$ is locally equivalent to an upper-block-triangular form iff  $d=d(E_1)+d(tr(E)_2)$ or $d=d(E_2)+d(tr(E)_1)$.
\\
\\
In the case of curves we have a much stronger criterion:

{\em Theorem \ref{Thm.Decomposability.Max.Gen.Curves.Upper.Block.Triang}
Let $(C,0)=(C_1,0)\cup(C_2,0)\sset(\k^2,0)$, with $(C_i,0)$ possibly further reducible, non-reduced but with no
common components. Let $\cM$ be a \dr\ of $(C,0)$.
\\1. $\cM\sim\bpm\cM_1&*\\\zero&\cM_2\epm$ iff  $\cM\sim\bpm\tilde{\cM}_2&*\\\zero&\tilde{\cM}_1\epm$,
where $\cM_i$, $\tilde{\cM}_i$ are some \dr s of $(C_i,0)$.
\\2. If $\cM$ is \mg\ at the origin then it is equivalent to an upper-block-triangular matrix, i.e. the corresponding module
is an extension.}
\\
\\
In general this extension of modules is non-trivial, the matrix is indecomposable. However the decomposability
holds if the components of the curve are not tangent.

{\em Theorem \ref{Thm.Decomposability.For.Tangential.Decomp.Curves}. Let $(C,0)=(C_1,0)\cup(C_2,0)\sset(\k^2,0)$
where $(C_i,0)$ can be further reducible or non-reduced but have no common tangents. If $\cM$ is a  \dr\ of this
curve that is \mg\ at the origin, then it is decomposable: $\cM\sim\cM_1\oplus\cM_2$.}

The last two results reduce the classification of local \mg\  \dr s of plane curve singularities
(i.e. families of matrices depending on two parameters) essentially to multiple branches, i.e. $(rC,0)$ with $(C,0)$
locally irreducible. In this case we have:

{\em Theorem \ref{Thm.Decomposability.Max.Gen.Det.Rep.Multiple.Curve}
Let $(rC,0)\sset(\k^2,0)$, where $(C,0)$ is a locally irreducible, reduced plane curve.
\\1. Let $\cM$ be a \dr\ of $(rC,0)$ \mg\ \underline{at the origin}. Then $\cM$ is equivalent to an upper-block-triangular matrix, the blocks
on the diagonal are \dr s of $(C,0)$.
\\2. Let $\cM$ be a \dr\ \mg\ \underline{on the punctured neighborhood} of the origin.
Then $\cM$ is totally decomposable: $\cM=\oplus\cM_i$ where $\cM_i$ is a \dr\ of $(C,0)$.}
\\
\\
In higher dimensional case we give an analog of theorem \ref{Thm.Decomposability.For.Tangential.Decomp.Curves}.
Let $(X,0)=(X_1,0)\cup(X_2,0)\sset(\k^n,0)$.

{\em Theorem \ref{Thm.Decomposability.For.Tangential.Decomp.Hypersurfaces}.
1. If the intersection $(X_1,0)\cap(X_2,0)$ is reduced, i.e. the components are reduced and
generically transverse, then any \dr\ that is \mg\ on the \underline{smooth points of $(X_1,0)\cap(X_2,0)$}
 is decomposable.
\\2. More generally, if the projectivized tangent cones,
$\P T_{(X_1,0)},\ \P T_{(X_2,0)}\sset\P(\k^n)$, intersect transversally then
any \dr\ of $(X,0)$ that is \mg\ \underline{near} \underline{the origin} is decomposable.}
\\\

As we show by numerous examples the assumptions of these criteria are almost necessary, so these sufficient
 criteria are in some sense the best possible.

\subsubsection{Relation to matrix factorizations}\label{Sec.Intro.Relations.To.MF}
Recall that a \mf\ of a function $f\in\cOn$ is the matrix identity $AB=f\one$. (We always assume $A|_0=\zero=B|_0$, for more detail see
\S\ref{Sec.Matrix.Factorizations.MCM}.) If $f$ is irreducible then $A$ is a \dr\ of some power of $f$.
 Which \dr s arise in this way? An immediate corollary of our approach is:

{\em Corollary \ref{Thm.Det.Reps.Matrix.Factorizat.Relation}.
  Let $\cM$  be a \dr\ of $\prod f_\al^{p_\al}$. It can be augmented to a \mf\ of $\prod f_\al$ (i.e.
  there exists $B$ with $MB=\prod f_\al\one$)
iff $\cM$ is \mg\ at smooth points of the reduced hypersurface $\{\prod f_\al=0\}\sset(\k^n,0)$}
\subsubsection{Saturated modules}\label{Sec.Intro.Contents.Saturated.Modules}
In \S\ref{Sec.Saturated.Det.Reps} we study \XXS \dr s. For a finite
modification $(X',0)\norm(X,0)$ the kernel $E$ of $\cM$ can be pulled back to $\nu^*E/Torsion$, a
module on $(X',0)$. Usually pulling back adds many new elements, as the initial kernel $E$ is not a
module over the bigger ring $\cO_{(X',0)}$. This gives a reformulation of definition \ref{Def.Saturated.Det.Reps}:

{\em Proposition \ref{Thm.Saturated.Modules.Properties}.
 The \dr\ $\cM$ is \XXS iff the corresponding kernel module $E$ is \XXS, i.e.
the canonical embedding $E\hookrightarrow\nu_*(\nu^*E/Torsion)$ is an isomorphism.
 The kernel $E$ of $\cM$ is \XXS iff the kernel $tr(E)$ of $\cM^T$ is \XXS.}

This proposition is helpful in proving various decomposability criteria.
For example:

{\em Theorem \ref{Thm.Sturated.Modules.Decomposability.Criterion}.
A \dr\ of $(X_1,0)\cup(X_2,0)$ is decomposable iff it is $\quotient{(X_1,0)\coprod(X_2,0)}{(X_1,0)\cup(X_2,0)}$
saturated.}
\\We emphasize that it is very simple to check saturatedness, by using
definition \ref{Def.Saturated.Det.Reps}, i.e. by checking the entries of $\cMv$.

Then we get an immediate:

{\em Corollary. Suppose $(X,0)$ is reduced and its normalization $(\tX,0)\norm(X,0)$ is a smooth variety.
(For example, the normalization of a reduced curve is always smooth.)
There exists unique $\tX/X$-saturated module (or $\tX/X$-saturated \dr): $\nu_*(\cO_{(\tX,0)})$. It is \mg.}
\subsection{A brief introduction and overview}\label{Sec.Intro Brief Overview}
We recall here the local aspects of \dr s only, for some references on the global aspects cf.\cite{Kerner-Vinnikov2009}.

\subsubsection{A view from singularities}\label{Sec.Intro View From Singularities}
The modern study started probably from the seminal paper \cite{Arnol'd1971} and
was mentioned in \cite[1975-26,pg.23]{Arnol'd-problems}.
Many works studied the miniversal deformations of a constant
matrix for various equivalences (i.e. to write a normal form for a linear family of matrices),
cf. e.g. \cite[Chapter5]{Tannenbaum81}, \cite{Khabbaz-Stengle70} or \cite{Lancaster-Rodman05}.
For the deformation theory from commutatie algebra point of view cf. \cite{Ile2004}.

Recently various singularity invariants of such "matrices of functions" have been established:
\cite{Bruce-Tari04}, \cite{Bruce-Goryun-Zakal02}, \cite{Goryun-Zakal03}, \\ \cite{Goryun-Mond05}.

\subsubsection{}\label{Sec.Intro View From Algebra} {\em The case of one variable} is elementary,
e.g. \cite[chapter VI]{Gantmacher-book}:
any square, non-degenerate matrix with formal entries is locally equivalent to the unique diagonal matrix:
$\oplus x^{d_i}\one_{r_i\times r_i}$ for $ d_1<\cdots<d_n$ and $\{r_i\}$ are some integers.
In more modern language: any $\k[[x]]$ module is the direct sum of cyclic modules.
\subsubsection{Matrix factorizations and maximally Cohen-Macaulay modules}\label{Sec.Matrix.Factorizations.MCM}
For an introduction to the case of more variables see \cite{Yoshino-book} and \cite{Leuschke-Wiegand-book}.

Let
 $\cM$ be a local \dr~ of  $f\in \cOn$. Let $E$ be its kernel
  spanned by the columns of $\cMv$ as a module over $\cOx:=\cOn/(f)$.
Then $E$ has a period two resolution by free $\cOx$ modules:
\beq
...\stackrel{\cMv}{\to}\cOx^{\oplus d}\stackrel{\cM}{\to}\cOx^{\oplus
d}\stackrel{\cMv}{\to}\cOx^{\oplus d}\to E\to0
\eeq
One can show that $depth(E)=n-1=dim\cOx$, hence $E$ is a maximally Cohen-Macaulay (MCM) module.
Maximally Cohen-Macaulay modules over reduced curves are just  the torsion-free modules.

Vice-versa \cite{Eisenbud80}: any maximally Cohen Macaulay (MCM) module $E$ over the hypersurface ring $\cOx$ as above
has a resolution of period two:
\beq
...\stackrel{\cM_1}{\to}\cOx^{\oplus d}\stackrel{\cM_2}{\to}\cOx^{\oplus d}\stackrel{\cM_1}{\to}
\cOx^{\oplus d}\to E\to0
\eeq
corresponding to the matrix factorization: $f\one=\cM_1\cM_2$.
Note that here the dimensions of $\{\cM_i\}$
are not fixed. So, in general $\cM_i$ are \dr s of $\prod f^{p_i}_i$, for $f_{reduced}=\prod f_i$.

Suppose $f$ is homogeneous, of degree $d$. Then, by \cite{BackHerzSand88}, $f$ admits a matrix factorization in
linear matrices: $f\one=A_1...A_d$, i.e. all the entries of $\{A_i\}$ are homogeneous linear forms.

For an MCM module $E$ over $\cOx$ the minimal number of generators of $E$ is not bigger
than $multiplicity(X,0)\times rank(E)$, \cite[\S3]{Ulrich84}.
Modules for which the equality occurs are called Ulrich's modules,
 In our case, with $rank(E)=1$, they are precisely the \mg\  \dr s.
For an arbitrary algebraic hypersurface Ulrich modules, of high rank, exist \cite[Theorem 1]{Backelin-Herzog89}.
Hence, for any $\{f=0\}\sset(\k^n,0)$ its multiples $\{f^p=0\}$, for $p$ high enough, have \mg\ \dr s.

Let $E$ be an MCM-module (of any rank) on $(X,0)=\{f=0\}$. Its resolution provides a \mf\ of $f$.
 The syzygy of $E$ gives a \dr\ of a "multiple" of $f$, i.e. a hypersurface $(Y,0)$ that contains $(X,0)$ and
 whose reduced locus coincides with the reduced locus of $(X,0)$. Thus $E$ can be considered as a rank one module
 on $(Y,0)$. So there exists a natural embedding of the theory of MCM-modules on hypersurfaces
 into the theory of rank one MCM-modules on (non-reduced) hypersurfaces.
\\
\\
MCM modules over a given hypersurface singularity, i.e. matrices of formal series with the given determinant,
 have been classified in some particular cases.
A hypersurface singularity is called of finite/tame CM-representation type if it has a finite/countable
number of indecomposable MCM's, up to isomorphism.

\li A series of papers resulted in \cite{Buchweitz-Greuel-Schreyer87}: a hypersurface singularity
 is of finite CM-representation type iff it is the ring of a simple (ADE) singularity.
\li The MCM's of rank 1 over locally irreducible plane curve singularities were thoroughly studied
in \cite{Greuel-Pfister-93}.
\li The MCM modules over the surface $\sum^3_{i=1} x^3_i$ and \mf s were classified in
\cite{LazaPfisterPopescu02}. The MCM modules of rank 2 over the hypersurface $\sum^4_{i=1} x^3_i$
 were classified in \cite{BaEnePfistPopesc05}.
\li The MCM modules over the ring $\quotient{k[[x,y]]}{(x^n)}$ were classified in \cite{Ene-Popescu08}.
\li The MCM modules over the ring $\quotient{k[[x,y]]}{(xy^2)}$  were classified in \cite{Buchweitz-Greuel-Schreyer87}
\li The MCM modules over surface singularities were studied in \cite{Burban-Drozd08}, in particular the modules
over $\quotient{\k[[x,y,z]]}{xy}$ and $\quotient{\k[[x,y,z]]}{x^2y-z^2}$ were classified. See also \cite{Burban-Drozd10}.
\li The MCM modules over Thom-Sebastiani rings, i.e. $\quotient{k[[x_1..x_k,y_1..y_n]]}{(f(x)\oplus g(y))}$,
were studied in \cite{Herzog-Popescu97}. In particular, the modules over $\quotient{k[[x_1,..,x_n,y]]}{(f(x)+y^3)}$ were
related to those over $\quotient{k[[x_1,..,x_n]]}{(f(x))}$. Some MCM modules over $\quotient{k[[x,y]]}{(x^a+y^3)}$
were classified.
\li The possible rank of an MCM module without free summand, on a reduced hypersurface is bounded from below.
In particular, the conjecture in \cite{Buchweitz-Greuel-Schreyer87} reads:
\beq
rank(E)\ge 2^{\lfloor\frac{\dim Sing(X,0)-2}{2}\rfloor}
\eeq
\subsubsection{Some applications}
\paragraph{} Semi-definite programming is probably the most important new development in optimization in the last two
decades. The general goal is to study the {\em linear matrix inequalities}, i.e. the positive-definiteness of
 $\sum x_i\cM_i$, where $\{\cM_i\}$ are real symmetric matrices.
 The literature on the subject is vast, for the general introduction cf. \cite{S.I.G.1997}.

A related object is the subsets of $\R^n$, presentable by linear matrix inequalities, called {\em spectrahedra}.
This notion was introduced and studied in \cite{RG95}, for the recent advances  cf.
\cite{Helton-Vinnikov2007}.

\paragraph{} The famous Lax conjecture relates the hyperbolic polynomials and self-adjoint positive
definite \dr s, \cite{Lax1958}. The initial form, for homogeneous polynomials in three variables has
been proved in \cite{Lewis-Parrilo-Ramana2005} and\\ \cite{Helton-Vinnikov2007}. The general case is still open.

\paragraph{} The spectral analysis of pairs of commuting non-selfadjoint
or non-unitary operators is essentially based on the determinantal representations of plane algebraic curves
cf. \cite{Livšic-Kravitsky-Markus-Vinnikov-book}.

\section{Preliminaries and background}\label{Sec.Background}
\subsection{On the base rings.}\label{Sec.Background.Base.Ring}
When studying \dr s several rings appear naturally:
\\- the ring of rational functions that are regular at the origin
(i.e. the localization at the the origin of the polynomial ring $\k[x_1..x_n]_{(\cm)}$),
\\- the ring of locally
 analytic functions $\k\{x_1..x_n\}$
\\-  the ring of formal power series $\k[[x_1..x_n]]$.

We always denote the maximal ideal by $\cm$.

The ring $\k\{x_1..x_n\}$ comes inevitably in the local considerations. For example, an algebraic hypersurface
singularity can be irreducible over $\k[x_1..x_n]_{(\cm)}$ but reducible over $\k\{x_1..x_n\}$.
The (ir)reducibility over $\k\{x_1..x_n\}$ and $\k[[x_1..x_n]]$ coincide, i.e. if $(X,0)$ is a locally
analytic germ and $(X,0)=\cup_i(X_i,0)$ is its decomposition into irreducible formal germs, then
all $(X_i,0)$ are locally analytic. This follows e.g. from Artin and Pfister-Popescu approximation theorems.

Even if one restricts to locally converging power series, in some inductive arguments the
formal power series might appear, fortunately just
as an intermediate step. The final result (the \dr\ and the matrices of equivalence)
can always be chosen locally analytic due to the approximation theorems:
\bthe\cite{Artin68}\cite[pg.32]{GLSbook}
Let $\underline{x},\underline{y}$ be the multi-variables and
 $f_1,..,f_k\in\k\{\underline{x},\underline{y}\}$ the locally analytic series.
 Suppose there exist formal power series
  $\hat{Y}_1(\underline{x}),..,\hat{Y}_l(\underline{x})\in\k[[\underline{x}]]$ solving the equations, i.e.:
\beq
f_i\Big(\underline{x},\hat{Y}_1(\underline{x}),..,\hat{Y}_l(\underline{x})\Big)\equiv0,~~~~i=1..k
\eeq
Then there exists a locally analytic solution ${Y}_1(\underline{x}),..,{Y}_l(\underline{x})\in\k\{\underline{x}\}$:
\beq
f_i\Big(\underline{x},{Y}_1(\underline{x}),..,{Y}_l(\underline{x})\Big)\equiv0,~~~~i=1..k
\eeq
\ethe
\bthe \cite{Pfister-Popescu75}
Let $F_1=0,F_2=0,..F_k=0$ be a system of polynomial equations over a complete local ring $(R,\cm)$.
The system has a solution in $R$ iff it has a solution in $R/\cm^N$ for any $N$.
\ethe

An immediate application of these theorems is:
\bcor
Let $\cM$ be a matrix with entries in $\k\{x_1,..,x_n\}$. Suppose for each $N$ the matrix,
 considered as matrix over $\k\{x_1,..,x_n\}/\cm^N$,
is equivalent to a block diagonal $\bpm *&\zero\\\zero&*\epm$ (or an upper-block-triangular $\bpm *&*\\\zero&*\epm$).
Then $\cM$ is equivalent to a block-diagonal (or an upper-block-triangular) matrix over $\k\{x_1,..,x_n\}$.
\ecor
Indeed, we have here the system of locally analytic equations, corresponding to the zero blocks of
 $A\cM B=\bpm *&\zero\\\zero&*\epm$ or  $A\cM B=\bpm *& *\\\zero&*\epm$. And these equations have a solution over
$\k\{x_1,..,x_n\}/\cm^N$ for any $N$.
\\
\\
Over the ring of locally analytic or formal power series the Krull-Schmidt theorem holds: the decomposition of $E$
into irreducible components (i.e. of $\cM$ into indecomposable blocks) is unique up to permutation.
For the situation over the rings of rational functions see \cite{Wiegand01} or \cite{Leuschke-Wiegand-book}.
\\
\\
Comparing the \dr s over various base rings we have two natural and well studied questions,
 cf. \cite{Eisenbud-book}, \cite{Yoshino-book}.
\li (injectivity) Let $\cM_1$, $\cM_2$ be \dr s with rational/locally converging entries.
Suppose they are formally equivalent. Are they rationally/locally converging equivalent?
This is indeed true, as the completion is a faithful functor.
\li (surjectivity) Which formal \dr s of an analytic/algebraic hypersurface are
equivalent to \dr s with analytic/rational entries?

By the classical theorem \cite[Th\'{e}or\`{e}me 3]{Elkik73} any formal \dr\ of an analytic hypersurface with
isolated singularity is equivalent to an analytic \dr.
For affine rings, i.e. the ring of rational functions, this property usually fails.

Any formal \dr\ of a zero dimensional hypersurface singularity, i.e. $f(x)\in\k[[x]]$,
 is equivalent to an algebraic (cf. the beginning of
\S\ref{Sec.Intro View From Algebra}).
For the \dr s of plane curve singularities one has:
\bpro
1. If $\det(\cM)$ defines an algebraic plane curve singularity that is locally irreducible (over $\k[[x,y]]$),
then $\cM$ is equivalent to a matrix of polynomials.
\\2. More generally, if $\det(\cM)$ defines an algebraic plane curve singularity all of whose
irreducible components (over $\k[[x,y]]$) are algebraic, then $\cM$ is equivalent to a matrix of polynomials.
\epro
This is proved e.g. in \cite[Proposition 3.3]{Frankild-Sather-Wagstaff-Wiegand08}. A
 down-to-earth proof is in \cite[Appendix A]{Belitskii-Kerner2010}

More generally, for any formal \dr\ $\cM$ of an algebraic hypersurface, there exists another, $\cN$,
such that $\cM\oplus \cN$ is equivalent to a \dr\ with polynomial entries,
see \cite[Corollary 3.5]{Frankild-Sather-Wagstaff-Wiegand08}.

\subsection{Hypersurface singularities}\label{Sec.Background.Hypersurface.Singularities}
\subsubsection{Some local decompositions}\label{Sec.Background.Decomposition into irreducibles}
For local considerations we always assume the (singular) point to be
at the origin and the ring is either $\k\{x_1,..,x_n\}$ or $\k[[x_1,..,x_n]]$.

Associated to any germ $(X,0)=\{f=0\}$ is the decomposition $(X,0)=\cup (p_iX_i,0)=\{\prod_i f^{p_i}_i=0\}$.
 Here each $(X_i,0)$ is reduced and locally irreducible. The {\em tangent cone} $T_{(X,0)}$ is formed
 as the limit of all the tangent planes at smooth points. Let $f=f_p+f_{p+1}+..$
be the Taylor expansion, then the tangent cone is $\{f_p=0\}\sset(k^n,0)$.
For curves the tangent cone is the collection of tangent lines,
each with the corresponding multiplicity.

The tangent cone is in general reducible, associated to it is the {\em tangential decomposition}:
$(X,0)=\cup (X_\al,0)$. Here $\al$ runs over all the (set-theoretical) components of the tangent cone,
each $(X_\al,0)$ can be further reducible, non-reduced.

The simplest invariant of the hypersurface singularity $\{f_p+f_{p+1}+\cdots=0\}$ is the {\em multiplicity} $p=mult(X,0)$,
 for the tangential components denote $p_\al=mult(X_\al,0)$.

\subsubsection{The normalization and intermediate modifications}\label{Sec.Background.Normalization.Intermediate.Modifications}
A singularity is normal if its local ring is a  domain integrally closed in its field of fractions,
\cite[pg.118]{Eisenbud-book}. A reduced hypersurface singularity
is normal iff its singular locus is of codimension at least two. For example, a reduced curve is normal iff it is smooth
and a reducible hypersurface is not normal.

The normalization of a (non-normal) germ is a (unique) finite proper birational morphism $(\tX,0)\to(X,0)$,
with $(\tX,0)$ pure dimensional and normal. Note that $(\tX,0)$ is usually a multi-germ, as the normalization
 separates the components. For brevity we write $(\tX,0)$ instead of $\coprod(\tX_i,0_i)$.
In the non-reduced case, $X=\cup (p_iX_i,0)$, the normalization is: $\coprod (p_i\tX_i,0)$.

Algebraically, if $\cOx$ is the local ring of $(X,0)$ then the normalization is induced by the inclusion
$\cOx\sset\cO_{(\tX,0)}$, where $\cO_{(\tX,0)}$ is the integral closure.

As any reduced normal curve is smooth, the normalization of a reduced curve singularity is its
resolution $(\tC,0)\to(C,0)$.
\\
\\
Usually the normalization $(\tX,0)\stackrel{\nu}{\to}(X,0)$ can be (nontrivially)
factorized: $(\tX,0)\to(X',0)\to(X,0)$. Here both maps are finite surjective bi-rational morphisms.
Usually this factorization can be done
in many distinct ways. All the possible intermediate steps form an oriented graph, usually not a tree.
 The initial vertex of this graph is the full normalization, the final is the original hypersurface.
Algebraically, the intermediate steps correspond to extensions of the local rings:
\beq
\cOx\stackrel{\nu_{\tX/X}^*}{\hookrightarrow}\cO_{(X',0)}\stackrel{\nu_{\tX/X'}^*}{\hookrightarrow}\cO_{(\tX,0)}
\eeq
\bex\label{Ex.A_n.Intermediate.Modifications} The curve singularity of type $A_n$: $y^2=x^{n+1}$.
\li $n=2l$. In this case the curve is a branch, i.e. is locally irreducible. The normalization is given by
$x=t^2$, $y=t^{2l+1}$ or by the extension of the local ring: $\k\{t^2,t^{2l+1}\}\subset\k\{t\}$. All the
 intermediate modifications correspond to the intermediate rings:
\beq
\k\{t^2,t^{2l+1}\}\sset\k\{t^2,t^{2l-1}\}\sset\k\{t^2,t^{2l-3}\}\sset..\sset\k\{t^2,t^3\}\sset\k\{t\}
\eeq
Or geometrically we have the surjections of plane curves singularities of types:
\beq
A_0\to A_2\to..\to A_{2l-2}\to A_{2l}
\eeq
$\bullet$ $n=2l-1$. In this case the curve has two branches. The normalization is given by
$x=t_1+t_2$, $y=t^l_1-t^l_2$ or by the extension of the local ring:
$\quotient{\k\{t_1+t_2,t^l_1-t^l_2\}}{\bl t_1t_2\br}\subset\k\{t_1\}\times\k\{t_2\}$. All the
 intermediate modifications correspond to the intermediate rings:
\beq
\quotient{\k\{t_1+t_2,t^l_1-t^l_2\}}{\bl t_1t_2\br}\sset\quotient{\k\{t_1+t_2,t^{l-1}_1-t^{l-1}_2\}}{\bl t_1t_2\br}
\sset...\sset\quotient{\k\{t_1,t_2\}}{\bl t_1t_2\br}\sset\k\{t_1\}\times\k\{t_2\}
\eeq
Or geometrically:
\beq
A_0\sqcup A_0\to A_1\to A_3\to..\to A_{2l-1}
\eeq
\eex
\bex\label{Ex.Intermediate.Modifications.for.D_4}
Consider the germ of the type $x^p=y^p$. This is the ordinary multiple point, i.e. the intersection of $p$ smooth
pairwise non-tangent branches. Here the tangent cone consists
of the lines $\{x=w y\}$, for $w^p=1$. The tangential decomposition coincides with the
branch decomposition. The normalization separates all the branches and is
defined by $(x=t_i,y=w_i t_i)$, here $w_i$ are all the $p$'th roots of unity.
This corresponds to the embedding $\quotient{\k\{x,y\}}{x^p=y^p}\sset\prod_i\k\{t_i\}$.

The graph of modifications for $p=3$ is:
\beq
A_0\sqcup A_0\sqcup A_0\to A_0\sqcup A_1\to Spec\Big(\quotient{\k\{t_1,t_2,t_3\}}{t_it_j,\ i\neq j}\Big)\to\{x^3=y^3\}
\eeq
Here $\quotient{\k\{t_1,t_2,t_3\}}{(t_it_j,\ i\neq j)}$ is the local ring of the curve singularity formed by three
pairwise non-tangent smooth branches. Its embedding dimension is 3, i.e. this is a non-planar singularity.

The term $A_0\sqcup A_1$ corresponds to the separation of one branch from the remaining two. (By permutation this
can be
done in three ways.) The term $A_0\sqcup A_0\sqcup A_0$ corresponds to the total separation of branches, i.e.
the normalization.
\eex
Note that in general most modifications lead to {\em non-planar} and even non-Gorenstein singularities.
\bex
A particular kind of modification is the separation of all the locally irreducible components: $\coprod(X_i,0_i)\to(\cup X_i,0)$.
It is isomorphism when restricted to each particular component.
For the rings: $\cO_{(\cup X_i,0)}\into\cO_{\coprod(X_i,0)}=\prod \cO_{(X_i,0_i)}$. If $E$ is a module
over $\cO_{(\cup X_i,0)}$ then it is lifted to the collection of
modules $\{E\!\!\underset{\cO_{(\cup X_i,0)}}{\otimes}\!\!\cO_{\coprod(X_i,0)}\}/Torsion$, defined by the diagonal
embedding $1_{(X,0)}\to\oplus 1_{(X_i,0)}$.
\eex

\subsubsection{Adjoint and conductor ideals}\label{Sec.Background.Adjoint.Conduct.Ideal}
Let $(X,0)\sset(\k^n,0)$ be a (possibly non-reduced) hypersurface singularity and $(X',0)\norm(X,0)$
a finite modification as above, e.g. the normalization. As $(X,0)$ is usually reducible, its modification
 is usually a multi-germ, $(X',0)=\coprod(X'_i,0_i)$, with the morphisms $(X'_i,0)\stackrel{\nu_i}{\to}(X,0)$.
\bed The relative conductor ideal is:
\beq
\cOx\supset I^{cd}_{X'/X}:=Ann_{\cOx}(\quotient{\cO_{(X',0)}}{\cOx})=\{g|~\forall i:~
\nu^*_i(g)\cO_{(X'_i,0)}\subset\nu^{-1}_i(\cOx)\}
\eeq
\eed
Note that by its definition $I^{cd}_{X'/X}$ is the maximal ideal both in $\cOx$ that is also an ideal in $\cO_{(X',0)}$.

Consider the ideals in $\cOn$ whose restriction to $(X,0)$
is contained in $I^{cd}_{X'/X}$. Call the maximal among them: the {\em relative adjoint ideal} $Adj_{X'/X}\sset\cOn$.
So, $Adj_{X'/X}|_{(X,0)}=I^{cd}_{X'/X}$. For various properties of the adjoint/conductor ideals for normalization in the reduced case
cf. \cite[\S IV.11]{Serre-book}, \cite{Fulton-02} \cite[I.3.4, pg 214]{GLSbook} and 

\bex\label{Ex.A_n.Conductor.Adjoint.Ideal}
Continue example \ref{Ex.A_n.Intermediate.Modifications}.
Let the local ring of $(C',0)$ be $\k\{t^2,t^{2l'+1}\}$, i.e. the modification is $A_{2l'}\to A_{2l}$. Then
\beq
I^{cd}_{C'/C}=\bl t^{2l-2l'},t^{2l+1}\br\sset\cO_{(C,0)}=\k\{t^2,t^{2l+1}\}, \ \ \
Adj_{C'/C}=\bl x^{l-l'},y\br\sset\cO_{(k^2,0)}
\eeq
Similarly, for the modification $A_{2l'-1}\to A_{2l-1}$ one has:
\beq
I^{cd}_{C'/C}=\bl t^{l-l'}_1+t^{l-l'}_2,t^l_1-t^l_2\br\sset\cO_{(C,0)}=\k\{t_1+t_2,t^l_1-t^l_2\}, \ \ \
Adj_{C'/C}=\bl x^{l-l'},y\br\sset\cO_{(k^2,0)}
\eeq
\eex
\bprop\label{Thm.Adjoint.Conductor.Ideals.for.Branch.Separation}
1. For the modifications $(X'',0)\to(X',0)\to(X,0)$ one has: $I^{cd}_{X''/X}\sset I^{cd}_{X'/X}$ and
$I^{cd}_{X''/X'}I^{cd}_{X'/X}\sset I^{cd}_{X''/X}\sset\cOx$.
\\2.
Let $(X,0)=(\cup X_i,0)$ where $(X_i,0)$ can be further reducible but with no common components.
Let $f$ and $\{f_i\}$ be the defining equations of $(X,0)$ and $\{(X_i,0)\}_i$.
Then \[Adj_{\coprod_i(X_i,0)/(\cup X_i,0)}=\bl \frac{f}{f_1},..,\frac{f}{f_k}\br\]
\\3. More generally, if $(X',0)=\coprod(X'_i,0_i)\norm(X,0)=(\cup X_i,0)$ is a finite modification then
\[I^{cd}_{X'/X}=\{\sum\frac{f}{f_i}g_i, \ g_i\in I^{cd}_{X'_i/X_i}\}\]
\eprop
\bpr 1. The product of any element of $I^{cd}_{X''/X'}$ with an element of $I^{cd}_{X'/X}$
lies in $\cOx$, hence one can speak about the inclusion.
The statement is immediate.
\\2. Suppose $g\in Adj_{\coprod_i(X_i,0)/(\cup X_i,0)}$,
so that $\nu^*(g)1_i\in\nu^{-1}\cOx$, for $1_i\in\cO_{(X_i,0)}$.  As $\nu^*(g)1_i|_{(X_j,0)}=0$ for $j\neq i$ we get that $\nu^*(g)1_i$
is divisible by  $\frac{f}{f_i}$. So $g\in\bl\frac{f}{f_i},f_i\br$. By going over all the components one get
$g\in\bl\frac{f}{f_1},..,\frac{f}{f_k}\br$.
So, $Adj_{\coprod_i(X_i,0)/(\cup X_i,0)}\sset\bl \frac{f}{f_1},..,\frac{f}{f_k}\br$. The converse inclusion is immediate.
\\3. The $\supset$ inclusion. Let $g=\sum\frac{f}{f_i}g_i$, with $\{g_i\in I^{cd}_{X'_i/X_i}\}_i$.
Let $h=\sum h_i\in\cO_{\coprod(X'_i,0_i)}$. Then $hg=\sum\frac{f}{f_i}g_ih_i$, as $h_jg_i=0$ for $j\neq i$.
Moreover $\{g_ih_i\in\cO_{(X_i,0)}\}_i$, hence $gh\in\cOx$.
\\
The $\sset$ inclusion. Let $g\in I^{cd}_{X'/X}$, so by the second part: $g=\sum\frac{f}{f_i}g_i$. Then
$\frac{f}{f_i}g_i\cO_{(X'_i,0)}\sset\cOx$. But for any $h_j\in\cO_{(X'_j,0)}$, with $j\neq i$:
$\frac{f}{f_i}g_ih_j|_{(X'_j,0)}=0$. Therefore $\frac{f}{f_i}g_ih_i\in\frac{f}{f_i}\cO_{(X,0)}$, because $\cOn$ is UFD.
Hence $g_ih_i\in\cOx$, i.e. $g_i\in I^{cd}_{X'_i/X_i}$.
\epr
\beR
$\bullet$ In general neither $I^{cd}_{C'/C}\sset\cO_{(C,0)}$ nor even $\nu^*I^{cd}_{C'/C}\sset\cO_{(C',0)}$
are principal ideals. For example for the modification
\beq
(C',0)=Spec(\k\{t^3,t^5,t^7\})\stackrel{\nu}{\to}(C,0)=Spec(\k\{t^3,t^7\}), \ \ \
\nu^*I^{cd}_{C'/C}=\bl t^7,t^9,..\br\sset\k\{t^3,t^5,t^7\}
\eeq
\li In general $I^{cd}_{C''/C'}I^{cd}_{C'/C}\subsetneq I^{cd}_{C''/C}$. For example,
\beq\ber
\cO_{(C,0)}=\k\{t^3,t^7\}\sset\cO_{(C',0)}=\k\{t^3,t^5,t^7\}\sset\cO_{(C'',0)}=\k\{t^3,t^4,t^5\}, \\
I^{cd}_{C'/C}=\bl t^7,t^9,..\br,\ I^{cd}_{C''/C'}=\bl t^3,t^5,..\br,\ I^{cd}_{C''/C}=\bl t^9,t^{10},..\br
\eer\eeq
\eeR

\subsection{The matrix and its adjoint}
We work with (square) matrices, their sub-blocks and particular entries. Sometimes to avoid confusion we
 emphasize the dimensionality, e.g. $\cM_{d\times d}$. Then $\cM_{i\times i}$ denotes an $i\times i$
block in $\cM_{d\times d}$ and $\det(\cM_{i\times i})$ the corresponding minor. By $\cM_{ij}$
we mean a particular entry.

Note that the adjoint of $\bpm \cM_1&\cM_3\\\zero&\cM_2\epm$ is
$\bpm\det(\cM_2)\cMv_1&-\cMv_1\cM_3\cMv_2\\\zero&\det(\cM_1)\cMv_2\epm$.

When working with matrices of functions, several natural notions arise:
\li $deg_{x_i}(\cM)$=the maximal degree of $x_i$ in the entries of $\cM$.
This is infinity unless all the entries of $\cM$ are
polynomials in $x_i$. Similarly for $deg(\cM)$, the total degree.
\li $ord_{x_i}(\cM)$=the minimal degree of $x_i$ appearing in $\cM$.
If an entry of $\cM$ does not depend on $x_i$ the order is zero, if $A\equiv\zero$ then $ord_{x_i}(A):=\infty$.
Similarly $ord(\cM)$  and $ord_x(\cM_{ij})$ for a particular entry.
So, e.g. $ord(\cM)\ge1$ iff $\cM|_0=\zero$
\li  $jet_k(\cM)$ is obtained from $\cM$ by truncation of all the monomials with total degree higher than $k$.

\subsubsection{Reduction to a minimal form}\label{Sec Local det reps Localization}
Let $\cM\in Mat(d\times d,\cOn)$, without the assumption $\cM|_{0}=\zero$.
Let the multiplicity of the hypersurface germ $\{\det(\cM)=0\}\subset(\k^n,0)$ be $mult(X,0)\ge1$.
\bpro\label{Thm.Localiz.Chip.off.Unity}
1. Locally $\cM_{d\times d}$ is equivalent to $\bpm \one_{(d-p)\times(d-p)}&\zero\\\zero& \cM_{p\times p}\epm$
with $\cM_{p\times p}|_{(0,0)}=\zero$ and $1\le p\le mult(X,0)$.
\\2. The stable equivalence (i.e. $\one\oplus\cM_1\sim\one\oplus\cM_2$) implies the ordinary local
equivalence ($\cM_1\sim\cM_2$).
\epro
This can be proved just by row and column operations of linear algebra.
From the algebraic point of view the first statement is the reduction to a minimal presentation of the module.
The second is the uniqueness of such a reduction. Both are proved e.g. in \cite[pg. 58]{Yoshino-book}.

The first statement is proved for the symmetric case in \cite[lemma 1.7]{Piontkowski2006}.
In the first statement both bounds are sharp, regardless of $d,p$ and $mult(X,0)$.
\subsubsection{Equivalence over $(\k^n,0)$ vs equivalence over $(X,0)$}
\bprop\label{Thm.Equivalence.Over.Hypersurface.Implies.Ordinary}
Two \dr s are equivalent over $(\k^n,0)$ iff they are equivalent over $(X,0)$.
\eprop
\bpr
The direct statement is trivial. For the converse statement, let $\cM_1\equiv A\cM_2B\mod\det(\cM_1)$.
One can assume both $\cM_1$ and $\cM_2$ vanish at the origin. Then $\cM_1=A\cM_2B+\cM_1\cMv_1Q$, for some
matrix $Q$ with entries in $\cOn$. Hence $\cM_1=A\cM_2B_2(\one-\cMv_1Q)^{-1}$. Here $(\one-\cMv_1Q)$ is
invertible as $\cMv_1|_0=\zero$.
\epr
Note that here $(X,0)$ is considered with its multiplicities, not just as a set. The equivalence
$\cM_1\sim\cM_2$ over $(X_{red},0)$ does not imply that over $(X,0)$. For example, both $\bpm x&0\\0&x^3\epm$
and $\bpm x^2&0\\0&x^2\epm$ restrict to zero matrix over $\{x=0\}$. But they are not equivalent when restricted
to $\{x^4=0\}$.
\subsubsection{The corank of the matrix}
Let $\cM$ be a determinantal representation of $(X,0)\sset(\k^n,0)$, let $\cMv$ be the adjoint matrix of $\cM$,
 so $\cM\cMv=\det(M)\one_{d\times d}$.

The corank of $\cM_{d\times d}$ at the point $pt\in\k^n$ is the maximal number $p$ such that
the determinant of any $(d-p+1)\times(d-p+1)$ minor of $\cM$ vanishes at $pt$.
The matrix $\cM$ is non-degenerate on $(\k^n,0)\smin(X,0)$ and the corank at the point $pt\in X$
 satisfies:
\beq\label{Eq Bound on Corank of Matrix}
1\le corank(\cM|_{pt})\le mult(X,pt)
\eeq
(To see this, note that $\cM$ is equivalent to $\one\oplus\cN$, where $\cN|_0=\zero$ and the corank of $\cM$ equals
the size of $\cN$.)
Hence any \dr~of a smooth hypersurface is \mg\ , cf. definition \ref{Def Maximality and Weak Maximality}.
For a reduced hypersurface the adjoint matrix $\cMv$ is not zero at smooth points of $X$.
As $\cMv|_X\times\cM |_X=\zero$ the rank of $\cMv$ at any smooth point of $X$ is 1.
 If $corank(\cM|_{pt})\ge2$ then $\cMv|_{pt}=\zero$.
Note that $\cM ^{\vee\vee}= (\det\cM)^{d-2}\cM $ and $\det\cMv=(\det\cM)^{d-1}$.

We have an immediate
\bprop If the representation is \mg\ at the origin, i.e. $rank(\cM_{d\times d}|_{0})=d-mult(X,0)$ then for
the minimal form $\cM_{p\times p}$: $p=mult(X,0)$
and  $det(jet_1\cM_{p\times p})\not\equiv0$ and $det(jet_{p-1}\cMv_{p\times p})\not\equiv0$.
\eprop

\subsubsection{Fitting ideals}
\bdp\label{thm Ideal Spanned by entries is invariant of local equivalence}
The fitting ideal $I_k(\cM)\sset\cOn$, is generated by all the $k\times k$ minors of $\cM$.
It is invariant under the local equivalence.
\edp
\bpr
First consider the case $k=1$, i.e. the ideal $I_1(\cM)$ is generated by the entries of $\cM$.
Then immediately: $I_1(A\cM B)\sset I_1(\cM)$. As $A,B$ are locally invertible the opposite inclusion holds too.

For arbitrary $k$ note that the wedge $\wedge^k\cM$ is the collection of all the $k\times k$ minors,
hence continue as for $k=1$.
\epr
\beR\label{Thm.Fitting.Ideal.Minimal.Number of Generators}
A trivial observation about the fitting ideals. Suppose $\cM_{p\times p}$ is locally
decomposable as $\bpm \cM_{p_1\times p_1}&\zero\\\zero&\cM_{p_2\times p_2}\epm$. Then
 $I_1(\cM)$ is generated by at most $p^2-2p_1p_2$ elements.

Similarly, if $\cM$ can be locally brought to an upper-block-triangular form then
$I_1(\cM)$ is generated by at most $p^2-p_1p_2$ elements.
\eeR
\ \\
\subsection{Kernel modules}\label{Sec.Background.Kernel.Modules}
Given a local \dr~ with $\cM|_0=\zero$, define the kernel module over $\cOn$ as follows.
Let $E\sset\cOn^{\oplus d}$ be the collection of all the kernel vectors, i.e. $v$ such that $\cM v$
is divisible by $det(M)$.
\bel
1. $E$ is a module over  $\cOn$, minimally generated by the columns of $\cMv$. It is supported on $(X,0)$.
\\2. Its restriction to the hypersurface (i.e. $E\otimes\cOx$) is a torsion-free module.
\\3. For the reduced hypersurface the module $E\otimes\cOx$ is free iff $\cM$ is a $1\times 1$ matrix.
\eel
\bpr
1. (This statement is also proved in \cite[pg.56]{Yoshino-book}.)
Let $E'$ be the  $\cOn$ module generated by the columns
of $\cMv$. Then $E'\sset E$. Let $v\in E$, so
$\cM v=\det(\cM)\bpm a_1\\..\\a_p\epm$. Let $v_1..v_p$ be the columns of $\cMv$,
then $\cM(v-\sum a_iv_i)=0\in\cOn$. As $\cM$ is non-degenerate on $(\k^n,0)$ we get $v\in E'$. Hence $E'=E$.
By linear independence on $(\k^n,0)$, the columns of $\cMv$ form a minimal set of generators.
\\2. The module $E_{(X,0)}$ is torsion-free as a submodule of the free module $\cO^{\oplus p}_{(X,0)}$.
\\3. Suppose $E$ is free and $v_1..v_p$ are the columns of $\cMv$, i.e. a minimal set of generators:
$E\approx \oplus\cOx v_i$.
Then $\cMv\cM|_{(X,0)}=\zero$ implies the linear relations among $\{v_i\}$, contradicting the freeness of $E$,
except in the case $\cM_{1\times 1}$.
\epr
By its definition the kernel module has a natural basis
$\{v_1..v_d\}$=the columns of $\cMv$. The embedded kernel with its basis determines the \dr:
\bpro\label{Thm.Kernels.Notions of equality}
1. Let $\cM_1,\cM_2\in Mat(d\times d,\cOn)$ be two local \dr~ of the same hypersurface
and $E_1,E_2$ the corresponding kernel modules.
Then $\cM_1=\cM_2$ or $\cM_1=A\cM_2$ or $\cM_1=A\cM_2B$ (for $A,B$ locally invertible) iff
$\Big(E_1,\{v^1_1..v^1_d\}\Big)=\Big(E_2,\{v^2_1..v^2_d\}\Big)\sset\cOn^{\oplus d}$
or $E_1=E_2\sset\cOn^{\oplus d}$
or $E_1\approx E_2$.
\\
\\2. In particular, if two kernel modules of the same hypersurface are abstractly isomorphic then
their isomorphism is induced by a unique ambient automorphism of $\cOn^{\oplus d}$.
\\3. In particular: $\cM$ is decomposable (or locally equivalent to an upper-block-triangular form) iff $E$
is a direct sum (or an extension).
\epro
Here in the first statement we mean the coincidence of the natural bases/the coincidence of the embedded
modules/the abstract isomorphism of modules.
\\\bpr {\bf 1,2.} The direction $\Rrightarrow$ in all the statements is immediate.
The converse follows from the uniqueness of minimal
free resolution \cite{Eisenbud-book}.

As the kernel is spanned by the columns of $\cMv$ the statement is straightforward,
except possibly for the last part: if $E_1\approx E_2$ then $\cM_1=A\cM_2B$.

Let $\phi: E_1\isom E_2$ be an abstract isomorphism of modules, i.e. an $\cOx$-linear map.
The module $E_1$ has a minimal free resolution. The isomorphism $\phi$ provides an additional minimal free resolution:
\beq\bM
0&\to &E_1&\to& \cO^{\oplus d}_{(X,0)}&\stackrel{\cM_1}{\to}&\cO^{\oplus d}_{(X,0)}...\\
&&\phi\downarrow&&\psi\downarrow&&\\
0&\to &E_2&\to& \cO^{\oplus d}_{(X,0)}&\stackrel{\cM_2}{\to}&\cO^{\oplus d}_{(X,0)}...
\eM\eeq
By the uniqueness of minimal free resolution, \cite[\S20]{Eisenbud-book}, we get the existence of
$\psi\in Aut(\cOx^{\oplus d})$.
\\
\\
{\bf 3.} Suppose $E=E_1\oplus E_2$, let $F_2\stackrel{\cM}{\to} F_1\to E\to 0$ be the minimal resolution.
Let $F^{(i)}_2\stackrel{\cM_i}{\to} F^{(i)}_1\to E_i\to 0$ be the minimal resolutions of $E_1,E_2$.
Consider their direct sum:
\beq
F^{(1)}_2\oplus F^{(2)}_2\stackrel{\cM_1\oplus\cM_2}{\to} F^{(1)}_1\oplus F^{(2)}_1\to E_1\oplus E_2=E\to 0
\eeq
This resolution of $E$ is minimal. Indeed, by the decomposability assumption the number of generators
of $E$ is the sum of
those of $E_1,E_2$, hence $rank(F_1)=rank(F^{(2)}_1)+rank(F^{(1)}_1)$.  Similarly, any linear relation
between the generators of $E$
(i.e. a syzygy)
 is the sum of relations for $E_1$ and $E_2$. Hence $rank(F_2)=rank(F^{(2)}_2)+rank(F^{(1)}_2)$.

Finally, by the uniqueness of the minimal resolution we get that the two proposed resolutions of $E$
are isomorphic, hence the statement.

Similarly for the extension of modules.
\epr
\subsubsection{Auslander transpose of $E$}\label{Sec.Background.Auslander.Transpose}
In addition to the kernel of $\cM$, spanned by the columns of $\cMv$, one sometimes considers the
left kernel: $tr(E)=Ker(\cM^T)$, called Auslander's transpose. It is spanned by the rows of $\cMv$.
 The propositions above hold for $tr(E)$
with obvious alterations. The modules $E$, $tr(E)$ are non-isomorphic in general. However (by the
last proposition) $E$ is decomposable or an extension iff $tr(E)$ is.
In addition, the minimal numbers of generators of $E$ and $tr(E)$ coincide.
\subsection{Pulback of modules: liftings and restrictions}\label{Sec.Background.Restriction.of.Module.To.Component}
Suppose a map of germs $(Y,0)\stackrel{j}{\to}(X,0)$ is given. In our context this will
be either a finite modification $(X',0)\twoheadrightarrow(X,0)$, as in the introduction, or the embedding $(X_1,0)\hookrightarrow(X_1\cup X_2,0)$,
here $(X_i,0)\sset(\k^n,0)$ are hypersurfaces, possibly reducible/non-reduced, but with no common components.

\subsubsection{Two ways to pullback}
The ordinary pull-back of a module is $E\underset{\cOx}{\otimes}\cO_{(Y,0)}$, usually it has torsion.
We always consider the torsion-free part: $E\underset{\cOx}{\otimes}\cO_{(Y,0)}/Torsion$.
\bex
\li Let $(C,0)=\{x^p=y^q\}\sset(\k^2,0)$ with $(p,q)=1$ and $q>p$. Consider the maximal ideal:
$\cm=\bl x,y\br\cO_{(C,0)}$.
Then the normalization $\nu:(\tC,0)\to(C,0)$ is $t\to(x=t^q,y=t^p)$ and $\nu^*(\cm)$ contains torsion. For
 example $\nu^*(x)-t^{q-p}\nu^*(y)$ is annihilated by $t^p=\nu^*(x)\in\k\{t\}$.
If we quotient by the torsion we get a free module: $\nu^*(\cm)/Torsion\approx\cO_{(\tC,0)}<\nu^*y>$.
\li Let $E$ be the kernel module of $\bpm f_1&0\\0&f_2\epm$, where $f_1,f_2\in\k\{x_1,..,x_n\}$ are mutually prime.
Then
\beq
E\underset{\quotient{\cOn}{(f_1f_2)}}{\otimes}\quotient{\cOn}{(f_1)}=\quotient{\cOn}{(f_1)}\bl s_1,s_2\br/(f_2s_2)
\eeq
while
\beq
\Big(E\underset{\quotient{\cOn}{(f_1f_2)}}{\otimes}\quotient{\cOn}{(f_1)}\Big)/Torsion=\quotient{\cOn}{(f_1)}\bl s_1\br
\eeq
\eex

In our case the kernel module is naturally embedded: $E_{(X,0)}\into\cOx^{\oplus d}$.
Hence the map $\cOx\stackrel{j^*}{\to}\cO_{(Y,0)}$ provides another version of pullback:
$j^*(i(E))\sset\cO^{\oplus d}_{(Y,0)}$. Here $j^*(i(E))$ is generated by the columns of $\cMv|_{(Y,0)}$.
The two pullbacks are compatible:
\bprop
The $\cO_{(Y,0)}$ modules $j^*(i(E))$ and $E\underset{\cOx}{\otimes}\cO_{(Y,0)}/Torsion$ are isomorphic.
\eprop
\bpr
Let $E_{(X,0)}$ be generated by $\{s_k\}_k$. Any element of $j^*(i(E))$ is presentable as $\sum a_k j^*(s_k)$,
where $a_k\in\cO_{(Y,0)}$. Thus a natural map is:
\beq
\phi:\ j^*(i(E))\ni\sum a_kj^*(s_k)\to[\sum a_k\otimes s_k]\in E\underset{\cOx}{\otimes}\cO_{(Y,0)}/Torsion
\eeq
This map is well defined. Indeed, if $\sum a_kj^*(s_k)=\sum b_kj^*(s_k)$ then
\beq
\phi(\sum a_kj^*(s_k))-\phi(\sum b_kj^*(s_k))=[\sum(a_k-b_k)\otimes s_k]=0\in E\underset{\cOx}{\otimes}\cO_{(Y,0)}/Torsion
\eeq
The map is linear and surjective by construction. Injectivity: if
$[\sum a_k\otimes s_k]=0\in E\underset{\cOx}{\otimes}\cO_{(Y,0)}/Torsion$, then there exists a non-zero divisor
$g\in\cO_{(Y,0)}$ such that $g(\sum a_k\otimes s_k)=0\in E\underset{\cOx}{\otimes}\cO_{(Y,0)}\sset\cO^{\oplus d}_{(Y,0)}$.
But then, as $g$ is a non-zero divisor, $\sum a_k\otimes s_k=0\in E\underset{\cOx}{\otimes}\cO_{(Y,0)}$.
 Hence the statement.
\epr

\subsubsection{Restriction to a component does not preserve Cohen-Macaulayness}
Suppose the hypersurface is locally reducible: $(X,0)=\cup(X_i,0)$. Consider the restriction
$E_i=E|_{(X_i,0)}/Torsion$, spanned by the column of $\cMv|_{(X_i,0)}$. By construction $E_i$ is
torsion-free, in particular if $(X,0)$ is a curve, i.e. $n=2$, then $E_i$ is Cohen-Macaulay.

In higher dimensions $E_i$ is not necessarily Cohen-Macaulay. For example
\beq
\cM=\bpm x&y\\0&z\epm\quad\quad\cMv=\bpm z&-y\\0&x\epm,\quad E=\quotient{\cOx\bl s_1,s_2\br}{(xs_1,ys_1+zs_2)}
\eeq
Here $E|_{x=0}$ is torsion free and isomorphic to the maximal ideal $\cm\sset\k\{x,y\}$. So it is a non-free module
over a regular ring, hence cannot be Cohen-Macaulay. On the other hand $E|_{z=0}/Torsion$ is free of rank 1.

In particular, in higher dimensions the minimal number of generators of $E_i$ can differ from that
of $tr(E)|_{(X_i,0)}/Torsion$.
\section{Decomposability of \mg\ \dr s}\label{Sec.Maximally.Generated.Reps.Decomposability.criteria}
Suppose the hypersurface is locally reducible $(X,0)=(X_1,0)\cup(X_2,0)$, where $(X_i,0)$ can be further
reducible, non-reduced but without common components.
Let $E$ and $E_i=E|_{(X_i,0)}/Torsion$ be the kernels of \dr s.
In this section we show that modules with large number of generators (e.g. \mg) tend to be extensions
or even decomposable.
\subsection{Preparations}\label{Sec.Maximally.Generated.Preparations}
\bel\label{Thm.Minor.Divisible.Det.is Very divisible} Let $\cM$ be an arbitrary square matrix with entries in $\cOn$.
\\1. Let $I\sset\cOn$ be a radical ideal that is a complete intersection. Suppose for any $i\times i$
minor of $\cM$ one has: $\det(\cM_{i\times i})\in I^l$. Then for any $(i+1)\times(i+1)$ minor:
$\det(\cM_{(i+1)\times(i+1)})\in I^{l+1}$.
\\2. In particular, suppose for any $i\times i$ minor $\cM_{i\times i}$ the determinant
is divisible by $g^l$ and $g$ has no multiple factors. Then, for any $(i+1)\times(i+1)$ minor
 $\cM_{(i+1)\times(i+1)}$, the determinant is divisible by $g^{l+1}$.
\\3. Consider the hypersurface germ $\{\prod_{i=1}^r f^{p_i}_i=0\}\sset(\k^n,0)$ for $\{f_i\}$ reduced.
Let $\cM$ be its \dr, $\cM|_{0}=\zero$. Suppose it is \mg\ on the locus $\cap_{j\in J}p_j X_j$
for $J\subseteq\{1,..,r\}$. Then all the entries of $\cMv$ belong to the power of the radical of the
ideal generated by $\{f_j\}_{j\in J}$: $\Big(Rad\bl \{f_j\}_{j\in J}\br\Big)^{\sum_{j\in J}p_j-1}$.
\\4. In particular, $\cM$ is \mg\ at the smooth points of the reduced locus $X_{red}$ iff all the entries
of $\cMv$ are divisible by $\prod f^{p_i-1}_i$.
\eel
\bpr
{\bf 1.} Let $\cM_{(i+1)\times(i+1)}$ be any minor, let $\cM_{(i+1)\times(i+1)}^\vee$ be its adjoint matrix.
By the assumption, every element of this adjoint matrix lies in $I^l$.
Hence
\beq
\Big(\det\cM_{(i+1)\times(i+1)}\Big)^i=\det\Big(\cM_{(i+1)\times(i+1)}^\vee\Big)\in I^{l(i+1)}
\eeq
Let $I=\bl g_1,..,g_k\br$ be a minimal set of generators. Consider the projection
$\cOn\to\quotient{\cOn}{\bl g_2,..,g_k\br}$. Let $\bl g_1\br\sset\quotient{\cOn}{\bl g_2,..,g_k\br}$ be the
image of $I$ under this projection. So, the image of $\Big(\det\cM_{(i+1)\times(i+1)}\Big)^i$ lies in
$\bl g_1^{l(i+1)}\br$. As $I$ is a complete intersection and $g_1$ is not a zero divisor one has:
\beq
\Big(\frac{\det\cM_{(i+1)\times(i+1)}}{ g_1^{l}}\Big)^i\in \bl g_1^{l}\br\sset\quotient{\cOn}{\bl g_2,..,g_k\br}
\eeq
As $g_1$ has no multiple factors one gets: the image of $\det\cM_{(i+1)\times(i+1)}$ in $\quotient{\cOn}{\bl g_2,..,g_k\br}$
 is divisible by $g_1^{l+1}$.

Finally note that the same holds for any generator of $I$. For example, for any $\k$-linear combination of the generators.
Hence the statement.
\\\\{\bf 2.} This is just the case of principal ideal, $I=\bl g\br$, for $g$ without multiple factors.
\\\\
{\bf 3.} Let $pt\in \cap_{j\in J}X_j$. By the assumption we have: $corank(\cM|_{pt})\ge\sum_{j\in J}p_j$.
So the determinant of any $(d-\sum_{j\in J}p_j+1)\times(d-\sum_{j\in J}p_j+1)$ minor of $\cM$ belongs to
the radical of the ideal generated by $\{f_j\}_{j\in J}$. By
the first statement we get: any entry of $\cMv$ belongs to $\Big(Rad\bl \{f_j\}_{j\in J}\br\Big)^{\sum_{j\in J}p_j-1}$.
\\\\
{\bf 4.} By the assumption, for any smooth point $pt\in\{f_i=0\}\smin\{\prod'_{j\neq i}f_j=0\}$
we have: $corank(\cM|_{pt})=p_i$. So, any
$(d-p_i+1)\times(d-p_i+1)$ minor of $\cM$ is divisible by $f_i$ near $pt$. By
the second statement we get: any entry of $\cMv$ is divisible by $f^{p_i-1}_i$ near $pt$.
Taking the closure we get the divisibility everywhere. Going over all the $\{f_i\}_i$ we get the direct statement.

For the converse statement, let $pt$ be a smooth point of the reduced locus. Can assume it is the origin.
Rectify the hypersurface
$\{\prod f_\al^{p_\al}=0\}$ locally
near this point, so the corresponding local ring is $\cOn/x_1^{p_\al}$. Restrict to the
line $x_2=0=..=x_n$. Then, analyzing the zero-dimensional case, one gets: the corank of $\cM$ at this
point is $p_\al$. Hence $\cM$ is \mg\ at the smooth points of the reduced locus.
\epr
\beR The conditions on the ideal in the proposition are relevant.
\li If the ideal is not a complete intersections the statement does not hold. For example, let $\cM_{3\times 3}$
be a matrix of indeterminates, let $I_2(\cM)$ be the ideal generated by all the $2\times 2$ minors.
One can check that $I_2(\cM)$ is radical. Hence any $2\times2$ minor belongs to $I_2(\cM)$ but
certainly $\det(\cM)\not\in I_2(\cM)^2$.
\li
In the third statement it is important to take \mg\  {\em near the point}. For example, $\cM=\bpm y&x\\0&y\epm$
is \mg\ at the origin but not near the origin. And not all the entries of $\cMv$ are divisible by $y$.
\eeR
Now we prove the statement of \S\ref{Sec.Intro.Relations.To.MF}
\bcor\label{Thm.Det.Reps.Matrix.Factorizat.Relation}
Let $\cM$ be a \dr\ of $\prod f_\al^{p_\al}$. It can be augmented to a \mf\ of $\prod f_\al$,
(i.e. there exists $B$ such that $\cM B=\prod f_\al\one$) iff $\cM$ is \mg\ at
the smooth points of the reduced locus $\{\prod f_\al=0\}$.
\ecor
\bpr $\Lleftarrow$ If $\cM$ is \mg\ at smooth points of the reduced locus then by theorem
\ref{Thm.Minor.Divisible.Det.is Very divisible} the adjoint matrix $\cMv$ is divisible by $\prod f_\al^{p_\al-1}$.
 Hence
\beq
\cM\frac{\cMv}{\prod f_\al^{p_\al-1}}=\prod f_\al\one
\eeq
$\Rrightarrow$ Suppose $\cM B=\prod f_\al\one$, for some matrix $B$.
Then $B=\frac{\prod f_\al}{\prod f_\al^{p_\al}}\cMv$, i.e. $\cMv$ is divisible by $\prod f_\al^{p_\al-1}$.
Now, by proposition \ref{Thm.Minor.Divisible.Det.is Very divisible}, $\cM$ is \mg\ at the smooth points of the
reduced locus.
\epr
\bthe\label{Thm.Decomposability.When.Adjoint.Matrix.Belongs.To.Separation.of.Branches}
Let $(X,0)=\{f_1f_2=0\}\sset(\k^n,0)$, where $f_i$ can be further reducible, non-reduced, but are relatively prime.
A \dr\ $\cM$ of $(X,0)$ decomposes as $\cM_1\oplus\cM_2$ iff every element of $\cMv$ belongs to the
ideal $\bl f_1,f_2\br\sset\cOn$.
\ethe
This proof uses only linear algebra. A more conceptual proof is in \S\ref{Sec.Saturated.Det.Reps}.
\bpr $\Rrightarrow$ is obvious.

$\Lleftarrow$
By the assumption $\cMv=f_2\cMv_1+f_1\cMv_2$, where $\cM_i$ are some (square) matrices with elements in $\cOn$,
Multiply the equality by $\cM$, then one has:
\beq
f_1f_2\one=\cM\cMv=f_2\cM\cMv_1+f_1\cM\cMv_2
\eeq
As $f_1$, $f_2$ are relatively prime, $\cM\cMv_i$ is divisible by $f_i$.
Therefore one can define the matrices $\{A_i\}$, $\{B_i\}$ by $f_iA_i:=\cM\cMv_i$ and $f_iB_i:=\cMv_i\cM$. By definition:
$A_1+A_2=\one$ and $B_1+B_2=\one$. We prove that in fact $A_1\oplus A_2=\one$ and $B_1\oplus B_2=\one$.
The key ingredient is the identity:
\beq
\cMv_j f_iA_i=\cMv_j\cM\cMv_i=f_jB_j\cMv_i
\eeq
Let $m_1,m_2$ be the multiplicities of $f_1,f_2$ at the origin.
It follows that $\cMv_jA_i$ is divisible by $f_j$ and thus $jet_{m_j-1}(\cMv_jA_i)=0$ for $i\neq j$. Hence, due to
the orders of $\cMv_j$,$\cM$ we get: $jet_{m_j}(\cM\cMv_jA_i)=jet_{m_j}(f_jA_jA_i)=0$, implying:
\beq
jet_0(A_j)jet_0(A_i)=jet_0(A_jA_i)\stackrel{for~i\neq j}{=}\zero,\text{  and  }\sum jet_0(A_i)=\one~~
\Rrightarrow ~~\one=\oplus jet_0(A_i)
\eeq
The equivalence $\cM\to U\cM X$ results in: $A_i\to U A_i U^{-1}$ and $B_j\to X B_j X^{-1}$.
So, by the conjugation by (constant) matrices can assume the block form:
\beq
jet_0(A_1)=\bpm \one&\zero\\\zero&\zero\epm,\ \ jet_0(A_2)=\bpm \zero&\zero\\\zero&\one\epm
\eeq
Apply further conjugation to remove the terms of $A_i$ in the columns of the i'th block to get:
\beq
A_1=\bpm \one&*\\\zero&*\epm,~~~~A_2=\bpm *&\zero\\ *&\one\epm
\eeq
Finally, use $A_1+A_2=\one$ to obtain $A_1=\bpm \one&\zero\\\zero&\zero\epm$ and $A_2=\bpm \zero&\zero\\\zero&\one\epm$.

Do the same procedure for $B_i$'s, this keeps $A_i$'s intact.
Now use the original definition, to write: $\cMv_i=\frac{f_i}{f}\cMv A_i$ and $\cMv_i=\frac{f_i}{f}B_i\cMv$.
This gives:
\beq
\cMv=f_2\cMv_1\oplus f_1\cMv_2
\eeq\epr
\subsection{Modules with many generators tend to be extensions}
Let $\cM_{d\times d}$ be a local \dr\ of a reducible hypersurface, $\{\det(\cM)=0\}=(X_1\cup X_2,0)$.
Here $(X_i,0)=\{f_i=0\}$ can be further reducible, non-reduced, but with no common components. As always, we
assume $\cM|_0=\zero$, i.e. the kernel $E$ is minimally generated by $d$ elements.

Suppose the restriction $E_i=E|_{(X_i,0)}/Torsion$, (cf.\S\ref{Sec.Background.Restriction.of.Module.To.Component}),
is minimally generated by $d(E_i)$ elements. In other words,
the maximal number of the columns of $\cMv|_{(X_i,0)}$, none of which belongs to the $\cO_{(X_i,0)}$
span of the others, is $d(E_i)$. Similarly, let $d(tr(E)_i)$ be the minimal number of generators for the
restrictions of the left kernels $tr(E)_i$, see \S\ref{Sec.Background.Auslander.Transpose}.
\bprop\label{Thm.Decomposability.Max.Gener.Can be Brought to Upper Triangular}
0. $d(E_i)>0$ and $\max(d(E_1),d(E_2))\le d\le d(E_1)+d(E_2)$. Similarly for $d(tr(E)_i)$.
\\1. $\cMv\sim\bpm f_1A_1&*\\\zero&f_2A_2\epm$, where $\zero$ is a $(d-d(tr(E)_2))\times(d-d(E_1))$ block of
zeros and $A_i$ some matrices with values in $\cOx$.
Similarly, $\cMv\sim\bpm f_2\tA_2&*\\\zero&f_1\tA_1\epm$, where $\zero$ is a $(d-d(tr(E)_1))\times(d-d(E_2))$
block of zeros.
\\2. In particular, $E$ is an extension (i.e. $\cM$ is equivalent to an upper block triangular) iff $d=d(tr(E)_1)+d(E_2)$
or $d=d(E_1)+d(tr(E)_2)$.
\eprop
\bpr
{\bf 0, 1.} The inequalities $0\le d(E_i)\le d$ are obvious. Consider $\cMv|_{(X_1,0)}$. By the assumption
the module of the columns of $\cMv|_{(X_1,0)}$ is generated by $d_1$ elements, hence one can assume
that the first $(d-d(E_1))$ columns of $\cMv$ are divisible by $f_1$. Similarly, the module of the rows of
$\cMv|_{(X_2,0)}$ is generated by $d(E_2)$ elements. Hence one can assume that the last $(d-d(E_2))$ rows
of $\cMv$ are divisible by $f_2$. As $f_1$, $f_2$ are relatively prime, one has:
\beq
\cMv\sim\bpm f_1A_1&*\\f_1f_2(..)&f_2A_2\epm
\eeq
Hence $\cMv|_{(X,0)}\sim\bpm f_1A_1&*\\\zero&f_2A_2\epm$. Now, by proposition
\ref{Thm.Equivalence.Over.Hypersurface.Implies.Ordinary},
 we have $\cMv\sim\bpm f_1A_1&*\\\zero&f_2A_2\epm$ over $(\k^n,0)$.

Similarly one obtains $\cMv\sim\bpm f_2\tA_2&*\\\zero&f_1\tA_1\epm$.

Finally, if $d_1+d_2<d$ then from the presentation above one gets: $\det(\cMv)\equiv0$ on $(\k^n,0)$.
\\\\
{\bf 2.} The direction $\Rrightarrow$ is obvious. For the converse statement,
by the first part we can assume $\cMv=\bpm f_1(..)&*\\\zero&f_2(..)\epm$, where $\zero$
is a $(d-d(tr(E)_2))\times(d-d(E_1))$ block of zeros. Hence if $d=d(E_1)+d(tr(E)_2)$
then $\cM=\bpm \cM_1&*\\\zero&\cM_2\epm$ with $\det(\cM_i)=f_i$.
\epr
\beR\label{Ex.Remark on the necessity of weak maximality}
$\bullet$ In general a matrix is not equivalent to an upper-block-triangular in two ways.
Namely, $\cM\sim\bpm\cM_1&*\\\zero&\cM_2\epm$, with $\cM_i$ a \dr\ of $(X_i,0)$, does not imply
$\cM\sim\bpm\tilde\cM_2&*\\\zero&\tilde\cM_1\epm$, with $\tilde\cM_i$ a \dr\ of $(X_i,0)$. For example
$\bpm x&y\\0&z\epm$ is not equivalent to $\bpm z&*\\0&x\epm$. For curves however this property does exist,
as is shown in theorem \ref{Thm.Decomposability.Max.Gen.Curves.Upper.Block.Triang}.
\li The assumption on the number of generators is essential. Consider the matrix
\beq
\cM=\bpm x^{p-1}y&x^p-y^p\\x^p+y^p&xy^{p-1}\epm,\quad p>2
\eeq
This \dr\ of an \omp~ is not locally equivalent to an upper triangular form.
Note that $I_1(\cM)$ is minimally generated by 4 elements, apply
remark \ref{Thm.Fitting.Ideal.Minimal.Number of Generators}.
\li
In general \mg\ \dr s are indecomposable. Consider
\beq
\cM=\bpm y+x^l&x^q\\0&y-x^l\epm,\quad 0<q<l
\eeq
To see that $\cM$ is not locally decomposable note that $I_1(\cM)=<y,x^q>$.
If $\cM\sim\bpm f_1&0\\0&f_2\epm$
then $f_i$ are the equations of branches of the curve $y^2=x^{2l}$ and $x^q\notin<f_1,f_2>$.
\li Even worse, being \mg\ at the origin does not imply equivalence to an upper-block-triangular. For example,
\beq
\cM=\bpm  x&-y&0\\z&0&y\\0&z&x\epm
\eeq
is a \dr\ of $\{xyz=0\}\sset(\k^3,0)$. The representation is \mg\ at the origin but is not equivalent to an
upper-block-triangular. If it were, the corank of $\cM$ would be at least $2$ on one of the intersections
$x=0=y$ or $x=0=z$ or $y=0=z$. But $corank(\cM)=1$ on all the intersections.
\eeR
\subsection{The case of curves}
For $n=2$ various strong criteria are possible.
\subsubsection{The criterion for being an extension}
\bthe\label{Thm.Decomposability.Max.Gen.Curves.Upper.Block.Triang}
Let $(C,0)=(C_1,0)\cup(C_2,0)\sset(\k^2,0)$, with $(C_i,0)$ possibly further reducible, non-reduced but with no
common components. Let $\cM$ be a \dr\ of $(C,0)$.
\\1. $\cM\sim\bpm\cM_1&*\\\zero&\cM_2\epm$ iff  $\cM\sim\bpm\tilde{\cM}_2&*\\\zero&\tilde{\cM}_1\epm$,
where $\cM_i$, $\tilde{\cM}_i$ are some \dr s of $(C_i,0)$.
\\2. If $\cM$ is \mg\ at the origin then it is equivalent to an upper-block-triangular matrix.
\ethe
\bpr
{\bf 1.} Let $E_i$ and $tr(E)_i$ be the restrictions of $Ker(\cM)$ and $Ker(\cM^T)$ to the local components, cf.
\S\ref{Sec.Background.Restriction.of.Module.To.Component}.
 Recall from \S\ref{Sec.Background.Auslander.Transpose} that for curves both $E_i$ and $tr(E)_i$ are Cohen-Macaulay,
i.e. torsion-free, hence their minimal number of generators coincide.

Suppose $\cM\sim\bpm \cM_1&*\\\zero&\cM_2\epm$, so $\cMv\sim\bpm f_2\cMv_1&*\\\zero&f_1\cMv_1\epm$. Here $\cM_i$ is
a $p_i\times p_i$ matrix.
We get: $tr(E)_1$ is minimally generated by $p_1$ elements. By the remark above: $E_1$ is minimally generated
by $p_1$ elements. Hence,the span of the columns of $\cMv|_{(X_1,0)}$ is generated by $p_1$ elements.
Thus $\cMv$ is equivalent to a matrix whose first $p_1$ columns are divisible by $f_1$. And this form can be
achieved by operations on columns only.

Similarly, $E_2$ is generated by $p_2$ elements, hence the same for $tr(E)_2$. Thus, after some row operations,
one can assume that the last $p_2$ rows on $\cMv$ are divisible by $f_2$.

Therefore $\cM_{(C,0)}$ has a $p_2\times p_1$ block of zeros, hence by property
\ref{Thm.Equivalence.Over.Hypersurface.Implies.Ordinary} the matrix $\cM$ has a block of zeros too.
\\\

{\bf 2.} Again, as the restrictions $E_i$ are Cohen-Macaulay, they are generated by at most $mult(C_i,0)$ elements.
Hence, the conditions of theorem \ref{Thm.Decomposability.Max.Gener.Can be Brought to Upper Triangular} are satisfied:
$d=mult(C,0)=d_1+d_2$ and the module is an extension.
\epr
\beR
In the first statement of the theorem, the matrices $\cM_i$ are in general not equivalent to $\tilde{\cM}_i$ or to
$\tilde{\cM}^T_i$. For example, consider a \dr\ of $y(y^2-x^{k+l})$:
\beq\ber
\bpm y&0&x\\0&y&x^l\\0&x^k&y\epm\rightsquigarrow\bpm y&0&x\\-x^{l-1}y&y&0\\0&x^k&y\epm\rightsquigarrow
\bpm y&0&x\\0&y&0\\x^{k+l-1}&x^k&y\epm\rightsquigarrow\bpm y&0&x\\x^{k+l-1}&y&x^k\\0&0&y\epm
\rightsquigarrow\bpm y&x&0\\x^{k+l-1}&y&x^k\\0&0&y\epm
\eer\eeq
And $\bpm y&x^l\\x^k&y\epm\not\sim\bpm y&x^{k+l-1}\\x&y\epm$.
\eeR
\subsubsection{Decomposability in the non-tangent case}
If the components of the curve are non-tangent then the \dr s tend to be decomposable.
 Recall the tangential decomposition $(C,0)=\cup_\al(C_\al,0)$
from \S\ref{Sec.Background.Decomposition into irreducibles}.
 Let $mult(C,0)=p$ and $mult(C_\al,0)=p_\al$. As always we assume $\cM|_0=\zero$.
\bthe\label{Thm.Decomposability.For.Tangential.Decomp.Curves}
Let $\cM_{p\times p}$ be a \dr\ of the plane curve $(C,0)$, \mg\ at the origin.
Corresponding to the tangential decomposition of $(C,0)$, the matrix $\cM$ is locally equivalent to a block diagonal:
$\cM_{p\times p}\sim\oplus_\al\cM_{p_\al\times p_\al}$. Here $\{\cM_{p_\al\times p_\al}\}$
are \mg\ \dr s of $\{(C_\al,0)\}$.
\ethe
\bpr
The theorem states that there exists a solution to the problem:
\beq
(\one+A)\cM(\one+B)=\bpm\cM_{m_1\times m_1}&\zero&..&\zero\\\zero&..\\\zero&..&\zero&\cM_{m_k\times m_k}\epm,~
A|_{(0,0)}=\zero=B|_{(0,0)}, ~ ~ \det(\cM_{m_\al\times m_\al})=f_\al
\eeq
for the unknowns $A,B,\{\cM_{m_\al\times m_\al}\}_\al$. Using Artin's and Pfister-Popescu theorems,
\S\ref{Sec.Intro View From Singularities}, it is enough to prove that the solution exists in
$\quotient{\k\{x,y\}}{\cm^N}$ for any $N$.

By the assumption $\cM $ vanishes at the origin, while the property \ref{Thm.Localiz.Chip.off.Unity} gives:
 $\det(jet_1\cM)\not\equiv0$.
\\{\bf Step 1.} By $GL(p,\k)\times GL(p,\k)$ bring  $jet_1(\cM)$ to the Jordan form.
For that, let $jet_1(\cM)=xP+yQ$ with $P,Q$ constant matrices. Assume that the curve is not tangent to coordinate axes.
Hence $P$ and $Q$ are of full rank.
By $GL(p,\k)\times GL(p,\k)$ bring $P$ to $\one$.
 The remaining transformation of $GL(p,\k)\times GL(p,\k)$
 preserving $P=\one$ is the conjugation: $\cM\to U\cM U^{-1}$. Hence $Q$ can be assumed in the Jordan form.
\\
\\
\\{\bf Step 2.}
The matrix $\cM $ is naturally subdivided into the blocks $B_{ij}$,
which are $p_i\times p_j$ rectangles (corresponding to the fixed eigenvalues of $jet_1(\cM)$).
We should remove the off-diagonal blocks, $B_{ij}$ for $i\neq j$.
We do this by induction, at the N'th step removing all the terms whose order is $\le N$.

Let $N=min_{i\neq j}(ord \cM_{ij})$ for $(ij)$ not in a diagonal block (thus $N>1$). Consider $jet_N(\cM )$,
i.e. truncate
all the monomials whose total degree is bigger than $N$. Suppose the block $B_{12}\subset jet_N(\cM )$ is non-zero,
i.e. there is an entry of order $N$.

As $l_1,l_2$ are linearly independent, by a linear change of coordinates in $(\k^2,0)$
can assume $l_1=x$, $l_2=y$. Decompose: $B_{12}=xT+yR$, where $T,R$
are $p_1\times p_2$ matrices, with $ord(T)\ge N-1$ and $ord(R)\ge N-1$.
From the last row of $B_{12}$ subtract  the rows
\beq
jet_N\cM_{p_1+1,*}, \ jet_N\cM_{p_1+2,*},..,jet_N\cM_{p_1+p_2,*}
\eeq
of $jet_N(\cM)$ multiplied
by $R_{p_11}$, $R_{p_12}$..$R_{p_1p_2}$. By the assumptions this does not change $jet_N(\cM)$ outside
the block $B_{12}$. After this procedure every entry of the last row of $B_{12}$
is divisible by $x$. Thus subtract from the columns of $B_{12}$ the column $jet_N\cM_{*,p_1}$ multiplied
by the appropriate factors.

Now the last row of $B_{12}$ consists of zeros, while $jet_N(\cM )$ is unchanged outside $B_{12}$.
Do the same procedure for the row $jet_N\cM_{p_1-1,*}$ of $B_{12}$ (using the
rows\\ $jet_N\cM_{p_1+1,*},jet_N\cM_{p_1+2,*},..,jet_N\cM_{p_1+p_2,*}$ and the column $jet_N\cM_{*,p_1-1}$). And so on.
\\
\\{\bf Step 3.}
After the last step one has the refined matrix $jet_N(\cM ')$ which coincides with $jet_N(\cM )$ outside
the block $B_{12}$ and has zeros inside this block. Do the same thing for all other (off-diagonal) blocks.
Then one has a block diagonal matrix  $jet_N(\cM ')$.

Now repeat all the computation starting from non-truncated version $\cM $. This results in the increase of $N$.
Continue by induction. Thus, for each $N$ can bring $\cM$ to such a form that the $jet_N(\cM )$ is
block diagonal. Then by the initial remark the statement follows.
\epr
\subsubsection{The case of multiple curve}
The results above reduce the decomposability questions to \dr s of a multiple
curve $(rC,0)\sset(\k^n,0)$, where $(C,0)$ is locally irreducible and reduced.
\bthe\label{Thm.Decomposability.Max.Gen.Det.Rep.Multiple.Curve}
Let $(rC,0)\sset(\k^2,0)$, where $(C,0)$ is a locally irreducible, reduced plane curve.
\\1. Let $\cM$ be its \dr\ \mg\ at the origin. Then $\cM$ is equivalent to an upper-block-triangular matrix, the blocks
on the diagonal are \dr s of $(C,0)$.
\\2. Let $\cM$ be a \dr\ \mg\ on the punctured neighborhood of the origin.
Then $\cM$ is totally decomposable: $\cM=\oplus\cM_i$ where $\cM_i$ is a \dr\ of $(C,0)$.
\ethe
\bpr Let $(\tC,0)\norm(C,0)=\{f=0\}$ be the normalization of the reduced curve. It defines valuation on $\cO_{(C,0)}$
by $val(g):=ord\nu^*(g)$, for $g\in\cO_{(C,0)}$.  In the non-reduced case the valuation on $\cO_{(rC,0)}$ is defined
by the pair:
\beq
val(g):=(ord_fg,val(\frac{g}{f^{ord_fg}}))
\eeq
Here $ord_f(g)$ is the maximal $k$ such that $g$ is divisible by $f^k$. In other words: $g\in f^k\cO_{(rC,0)}$
but $g\not\in f^{k+1}\cO_{(rC,0)}$.

The natural order on pairs for this valuation is defined by
\beq
(a_1,a_2)<(b_1,b_2)\quad if\ a_1<b_1\quad or\ \ \bpm a_1=b_1\\a_2<b_2\epm
\eeq
Let $p=mult(C,0)$, so $\cM$ is a $pr\times pr$ matrix.

{\bf 1.} Compare the valuations of the entries of $\cMv$. After a permutation of rows and columns we can assume
that $\cMv_{1,1}$ has the minimal valuation in $\cMv$ and in the first row and column the valuations are increasing:
\beq
val(\cMv_{1,1})<val(\cMv_{1,2})<..<val(\cMv_{1,pr}),\quad\quad
val(\cMv_{1,1})<val(\cMv_{2,1})<..<val(\cMv_{pr,1})
\eeq
Note that $\cMv_{11}\not=0\in\cO_{(rC,0)}$. Note that $\cO_{(\tC,0)}$ is generated as a $\cO_{(C,0)}$
module by $p$ elements. Therefore we can assume (possibly after a subtraction of columns) that the
elements $\cMv_{1,j}$ are divisible by $f^{\lfloor\frac{j}{p}\rfloor}$. Similarly, after some row subtraction
we can assume that the elements $\cMv_{j,1}$ are divisible by $f^{\lfloor\frac{j}{p}\rfloor}$.

Finally, recall that $rank(\cMv|_{(rC,0)})\le1$, i.e. any two rows or columns are proportional.
Hence, in the chosen basis the matrix is:
\beq
\cMv=\bpm *&*&*&*\\ * &*&*&\zero\\..&..&..&..\\ * &\zero&..&\zero\epm
\eeq
i.e. is equivalent to an upper-block-triangular.

{\bf 2.} By Proposition \ref{Thm.Minor.Divisible.Det.is Very divisible} the adjoint matrix $\cMv$ is divisible by $f^{r-1}$.
 Let $\cNv_{p\times pr}$ be the submatrix of $\frac{\cMv}{f^{r-1}}$ formed by lower $p$ rows.
Consider the module over $\cO_{(C,0)}$ spanned by the columns of $\cNv$. This module is generated by $p$ elements.
This can be seen, for example, by checking the valuation of the columns, by $\tC\norm(C,0)$.

Hence the matrix $\frac{\cMv}{f^{r-1}}$ is equivalent to the upper-block-triangular matrix, with the
zero block $\zero_{p\times(r-1)p}$. Hence $\cMv$ is equivalent to the upper-block-triangular matrix. Assume $\cMv$
in this form. Now consider the submatrix of $\frac{\cMv}{f^{r-1}}$ formed by the last $p$ columns.
By the argument as above one gets: $\cMv$ is equivalent to a block diagonal, with blocks $p\times p$
and $(r-1)p\times(r-1)p$.

Continue in the same way to get the statement.
\epr

\subsection{Higher dimensional case}\label{Sec.Decomposability.Higher.Dim.Case}
In some cases we have decomposability according to the tangential decomposition of a reduced hypersurface.
\bthe\label{Thm.Decomposability.For.Tangential.Decomp.Hypersurfaces}
Let $n\ge3$ and $(X,0)=(X_1,0)\cup(X_2,0)$.
\\1. If the intersection $(X_1,0)\cap(X_2,0)$ is reduced, i.e. the components are reduced and
generically transverse, then any \dr\ that is \mg\ on the smooth points of $(X_1,0)\cap(X_2,0)$
 is decomposable.
\\2. More generally, if the projectivized tangent cones,
$\P T_{(X_1,0)},\ \P T_{(X_1,0)}\sset\P(\k^n)$, intersect transversally then
any \dr\ of $(X,0)$ that is \mg\ near the origin is decomposable.
\ethe
\bpr {\bf 1.} By part 3 of proposition \ref{Thm.Minor.Divisible.Det.is Very divisible} every entry of $\cMv$
belongs to $\bl f_1,f_2\br\sset\cOn$. Then the decomposability follows by proposition
 \ref{Thm.Decomposability.When.Adjoint.Matrix.Belongs.To.Separation.of.Branches}.
\\{\bf 2.} As $\cM$ is maximally generated near the origin, for any point $pt\in X_1\cap X_2$
the order of vanishing of any element of the adjoint matrix satisfies: $ord_{pt}\cMv_{ij}\ge mult(X,pt)-1$.
 We claim that this implies $\cMv_{ij}\in\bl f_1,f_2\br\sset\cOn$, hence as above $\cM$ is decomposable
 (by proposition \ref{Thm.Decomposability.When.Adjoint.Matrix.Belongs.To.Separation.of.Branches}).

So, we should prove the following statement: given $f_1,f_2,h\in\cOn$, such that the projectivized
tangent cones $\P T_{(f_1=0)},\ \P T_{(f_2=0)}\sset\P(\k^n)$ intersect transversely and for any
point $pt\in(\k^n,0)$: $ord_{pt}(h)\ge ord_{pt}(f_1)+ord_{pt}(f_2)-1$. Then $h\in\bl f_1,f_2\br\sset\cOn$.

Let $Z=\{f_1=0=f_2\}$, then at any point of $Z$ the order of $h$ is at least one, i.e. $h$ vanishes on $Z$,
i.e. $h\in Rad\bl f_1,f_2\br\sset\cOn$. We should prove that $h$ belongs to the ideal $\bl f_1,f_2\br$ itself.
The proof is by induction on $n$.

Suppose $n=2$, i.e. $\{f_i=0\}$ are curve singularities whose tangent cones intersect at the origin only.
Let $H^0(\cO_{(\k^2,0)}(d))$ be the vector space of all the homogeneous polynomials in two variables of degree $d$.
Let $H^0(\cO_{(\k^2,0)}(-f_i)(d))\sset H^0(\cO_{(\k^2,0)}(d))$ be the subspace of all the polynomials divisible
by $f_i$. Then we have the exact sequence:
\beq\small
0\to H^0(\cO_{(\k^2,0)}(-f_1-f_2)(d))\to H^0(\cO_{(\k^2,0)}(-f_1)(d))\oplus H^0(\cO_{(\k^2,0)}(-f_2)(d))\to
H^0(\cO_{(\k^2,0)}(d))
\eeq
By the assumption $ord_0h\ge ord_0(f_1)+ord_0(f_2)-1$.
Hence it is enough to show that the map $H^0(\cO_{(\k^2,0)}(-f_1)(d))\oplus H^0(\cO_{(\k^2,0)}(-f_2)(d))\to
H^0(\cO_{(\k^2,0)}(d))$ is surjective for $d\ge  ord_0(f_1)+ord_0(f_2)-1$.
This is checked by computing the dimensions:
\beq\ber\small
dim H^0(\cO_{(\k^2,0)}(d))=d+1,\quad\quad
dim H^0(\cO_{(\k^2,0)}(-f_i)(d))=d-d_i+1,\\
dim H^0(\cO_{(\k^2,0)}(-f_1-f_2)(d))=d-d_1-d_2+1
\eer\eeq
here $d_i=ord(f_)$.
Hence for $n=2$ we get: $h\in\bl f_1,f_2\br\sset\cOn$.

Suppose the statement has been proven for the case of $(n-1)$ variables. Let $pt\in Z=X_1\cap X_2$, let $(\k^{n-1},pt)\sset(\k^n,pt)$
be a hyperplane transversal to the tangent cones $T_{(X_1,pt)}$, $T_{(X_2,pt)}$. Then the
restrictions $f_1|_{(\k^{n-1},pt)}$, $f_2|_{(\k^{n-1},pt)}$, $h|_{(\k^{n-1},pt)}$ satisfy the assumptions
of the statement. Hence by the induction assumption: $h=a_1f_1+a_2f_2+lh'$, where $a_1,a_2,h'$ are
some regular functions, while $l$ is the locally defining equation of the hyperplane $(\k^{n-1},pt)$.
Note that $h'$ itself satisfies the assumptions of the statement on the punctured neighborhood of $pt\in Z$.
As the vanishing order does not increase under small deformations we get that $h'$ satisfies the assumptions
of the statement on the whole neighborhood of $pt\in Z$. Then reiterating procedure we get
$h=a'_1f_1+a'_2f_2+l^2h''$ etc. As $h$ is locally analytic, this process stops after a finite number of steps,
giving $h\in\bl f_1,f_2\br\sset\cOn$.
\epr
\beR Note that for curves (theorem \ref{Thm.Decomposability.For.Tangential.Decomp.Curves}) we ask for \mg\
{\em at the origin},
while in higher dimensions (theorem \ref{Thm.Decomposability.For.Tangential.Decomp.Hypersurfaces}) we ask
for \mg\ on an open set {\em near the origin}. This is essential. For example $\cM=\bpm x&y\\0&z\epm$ is \mg\
at the origin.
And the hypersurface $\{\det(\cM)=xz=0\}$ consists of two transverse hyperplanes.
 But the \dr\ is not \mg\ {\em near} the singular point and is indecomposable. A similar example is in remark
 \ref{Ex.Remark on the necessity of weak maximality}.
\eeR
The theorem implies an immediate
\bcor
Let $(X,0)=\cup_\al(X_\al,0)$ be the reduced union of pairwise non-tangent smooth hypersurfaces,
e.g. an arrangement of hyperplanes.
Then $(X,0)$ has the unique \dr\ \mg\ on the neighborhood of $0\in\k^n$: the diagonal matrix.
\ecor

\subsection{Limits of kernel fibres}
One often imposes the following conditions of linear independence. Let $(X,0)=\cup(X_i,0)$ be reduced,
$E_X\sset X\times\k^d$ and
$E_i=E|_{(X_i,0)}/Torsion$. Let $Y_i\sset X_i$ be the maximal subvariety over which $E_i$ is locally free.
(So $Y_i$ is open dense in $X_i$.) Consider the topological closure of the embedded line bundle:
$\overline{E_{Y_i}}\sset X_i\times\k^d$. Denote its fibre at the origin by $\overline{E_{Y_i}}|_0$.

\bprop\label{Thm.Localiz.Final.Decomposability}
Let $(X,0)=\cup(X_i,0)\sset(\k^n,0)$ be a collection of reduced, smooth hypersurfaces.
The \dr\ $\cM_{p\times p}$ is completely decomposable iff the fibres $\{\overline{E_{Y_i}}|_0\}$ are
one dimensional vector subspaces of $\k^d$ that are linearly independent:
$Span(\cup \overline{E_{Y_i}}|_0)=\oplus \overline{E_{Y_i}}|_0$.
\eprop
\bpr
$\Rightarrow$ is obvious.
\\$\Leftarrow$  By the assumption $\cM|_{0}=\zero$ hence the corank of $\cM|_{0}$ equals the number of
(smooth) components, i.e. the multiplicity of the singularity. By continuity of the fibres (embedded vector spaces)
this happens also at the neighboring points. Hence $\cM$ is \mg\ near the origin.

Now, fix $\overline{E_{Y_1}}|_0=(1,0,..,0)\in\k^d$ and "rectify" the fibres locally. Namely, after a $GL(\k^d)$
transformation one can assume: $\overline{E_{Y_1}}=(1,0,..,0)$ over some neighborhood of the origin, while
$\overline{E_{Y_{i>1}}}\sset\{z_1=0\}\sset\k^d$ near the origin.
Hence, in this basis
\beq
\cMv=\bpm *& \zero&..&\zero\\\zero&* &*&*\\..&..&..&..\\\zero&*&*&*\epm
\eeq
Repeat for other components.
\epr
\beR\label{Ex.Kernels.Lin.Indep.NonSmooth.Branches.Indecomposable}
It is not clear whether the conditions can be weakened.
\li The fibre at the origin, $\overline{E_{Y_i}}|_0$, can be not a one-dimensional vector space, cf. the last
example in remark \ref{Ex.Remark on the necessity of weak maximality}.
\li The smoothness of the components in the
statement is important.
For example, consider $\cM=\bpm x^a& y^{d+1}\\y^c& x^by\epm$, a \dr~ of $y(x^{a+b}-y^{c+d})$ for $(a+b,c+d)=1$.
Assume also $c>1$ and $d>0$.
This \dr~ is not equivalent to an upper-triangular. Otherwise one would have $I_1(\cM)\ni y$.

On the other hand the limits of the kernel sections are linearly independent.
$\cMv=\bpm x^by&-y^{d+1}\\-y^c&x^a\epm$. So, on $y=0$ the kernel
is generated by $\bpm 0\\x^a\epm$, whose limit is $\bpm 0\\1\epm$. On $x^{a+b}=y^{c+d}$
both columns of $\cMv$ are non-zero, but linearly dependent.  So, for $a> d+1$ or $c-1>b$ their (normalized)
limit at the origin is $\bpm 1\\0\epm$.
\li It is important to ask for the {\em common} linear independence of the fibres, not just the pairwise
linear independence. Recall the first example in remark \ref{Ex.Remark on the necessity of weak maximality}.
There the branches are smooth and the limits of any two fibres are independent.
But altogether they are not linearly independent.
\eeR

%
%
%
%

%
%
%
%
%

%
%
%
%
%
\section{Saturated \dr s}\label{Sec.Saturated.Det.Reps}
Let  $(X',0)\norm(X,0)$ be a finite modification (cf. introduction).
Given a torsion-free module $E_{(X,0)}$, let $\nu^*E_{(X,0)}/Torsion$ be its torsion free pull-back
(\S\ref{Sec.Background.Restriction.of.Module.To.Component}).
 Then $E$ is naturally embedded into the pushforward $E\sset\nu_*(\nu^*E/Torsion)$.
\bed\label{Def.Saturated.Modules}
The module $E_{(X,0)}$ is called $X'/X$ saturated if this embedding is an isomorphism of modules:
$E\isom\nu_*(\nu^*E/Torsion)$
\eed
\bex\label{Ex.A_n.Minimal.Liftings}
Consider the torsion-free modules of rank 1 over the $A_n$ singularity $y^2=x^{n+1}$, continuing
examples \ref{Ex.A_n.Intermediate.Modifications} and \ref{Ex.A_n.Conductor.Adjoint.Ideal}.
Every such module can be embedded as $\cO_{(C,0)}\subset E\subset\cO_{(\tC,0)}$.
\li $n=2k$. The torsion-free modules of rank 1 are $\big\{E_l:=\cO_{(C,0)}\bl 1,t^{2l+1}\br\big\}_{0\le l<k}$.
Let $\cO_{(C',0)}=\k\{t^2,t^{2l+1}\}$ and $(C',0)\norm(C,0)$ the corresponding modification.
Then $\nu^*(E_l)/Torsion$ is a free $\cO_{(C',0)}$ module. In fact $E_l=\nu_*\cO_{(C',0)}$.
The corresponding \dr\ is $\bpm y& x^{k-l}\\ x^{l+k+1}& y\epm$.
\li $n=2k-1$. The torsion-free modules of rank 1 are $\{E_l:=\cO_{(C,0)}\bl 1,t^l_1-t^l_2\br\}_{0\le l<k}$.
Let $\cO_{(C',0)}=\k\{t_1+t_2,t^l_1-t^l_2\}$ and $(C',0)\norm(C,0)$ the corresponding modification.
Then $E_l=\nu_*\cO_{(C',0)}$. The corresponding \dr\ is $\bpm y+x^k& x^{k-l}\\0& y-x^k\epm$.

Note that in both cases the ideal generated by the entries of $\cMv$ is precisely the relative
adjoint ideal $Adj_{C'/C}$.
\eex

\bprop\label{Thm.Saturated.Modules.Properties} 0. Every torsion-free module over $(X,0)$ is $X/X$ saturated.
\\1. For any torsion-free module $E_{(X,0)}$ there exists the unique maximal finite modification
such that $E$ is \XXS. Namely, if $E$ is also $X''/X$ saturated then the modification $(X',0)\to(X,0)$
factorizes as $(X',0)\to(X'',0)\to(X,0)$.
\\2. Let $E$ be the kernel of a \dr\ $\cM$ of $(X,0)$. Then $E$ is \XXS iff $\cM$ is \XXS,
in the sense of definition \ref{Def.Saturated.Det.Reps}.
\\3. $E$ is \XXS iff $tr(E)$ is \XXS.
\eprop
\bpr {\bf 0.} Trivial.\\
\parbox{12.9cm}
{{\bf 1.} Suppose $E$ is both $X_1/X$ and $X_2/X$ saturated for the extensions of local rings:
$\cOx\sset\cO_{(X_1,0)},\cO_{(X_2,0)}\sset\cO_{(\tX,0)}$, here $\tX$ is the normalization.
Let $R\sset\cO_{(\tX,0)}$ be the subring generated by $\cO_{(X_1,0)},\cO_{(X_2,0)}$. Geometrically we
have the diagram on the right.}\hspace{0.5cm}
$\bM Spec(R)&\to&(X_2,0)\\\da&&\da\\(X_1,0)&\to&(X,0)\eM$
\\\\
 Now, by construction, $E$ is $Spec(R)/X$ saturated. By taking such unions of the
local rings (and staying inside $\cO_{(\tX,0)}$) the unique maximal modification is constructed.
\\
\\{\bf 2.} $\Rightarrow$ If $E$ is \XXS then it is an $\cO_{(X',0)}$ module. Recall that $E$ is spanned
by the columns of $\cMv$. Hence for any entry of $\cMv$: $\cO_{(X',0)}\cMv_{ij}\in\cOx$, i.e. $\cMv_{ij}\in Adj_{X'/X}$.
\\ $\Leftarrow$ If all the entries of $\cMv$ belong to $Adj_{X'/X}$ then $\nu_*\nu^*(E)/Torsion$ is generated
(as an $\cOx$ module) by some columns with entries in $\cOx$. Let $s\in\nu_*\nu^*(E)/Torsion$,
then $\cM s=0\in\cOx^{\oplus d}$. But then, by definition, $s\in E$. Hence the statement.
\\{\bf 3.} Note that $\cM$ is \XXS iff $\cM^T$ is.
\epr
The notion of being $X'/X$ saturated suits perfectly for the decomposition criterion.
\bthe\label{Thm.Sturated.Modules.Decomposability.Criterion}
Let $(X,0)=(X_1,0)\cup(X_2,0)\sset(\k^n,0)$ where $(X_i,0)=\{f_i=0\}$ can be further reducible, non-reduced, but without
common components. Let $X'=(X_1,0)\coprod(X_2,0)\to(X,0)$ be the finite modification that separates the components.
Let $E$ be the kernel of a \dr\ $\cM$ of $(X,0)$. The following are equivalent.
\\1. $E$ is $X'/X$ saturated.
\\1'. $\cM$ is $X'/X$ saturated, i.e. every element of $\cMv$ belongs to the ideal $\bl f_1,f_2\br\sset\cOn$.
\\2. $E=E_1\oplus E_2$, where $E_i=E|_{(X_i,0)}/Torsion$.
\\2'. $\cM\sim\cM_1\oplus\cM_2$, where $\cM_i$ is a \dr\ of $(X_i,0)$.
\ethe
In particular, if $(X,0)=\{f=0\}$, with $f=\prod f_i$, and $\cMv=\sum\frac{f}{f_i}\cMv_i$
 then the \dr\ is completely decomposable: $\cM\sim\oplus\cM_i$.

 Note that (for the statement 1') the relative adjoint ideal of $Adj_{\quotient{(X',0)}{(X,0)}}$
was computed in Proposition \ref{Thm.Adjoint.Conductor.Ideals.for.Branch.Separation}.
\\\bpr
The equivalence of 1 and 1' is proven in Proposition \ref{Thm.Saturated.Modules.Properties}.
The equivalence of 2 and 2' is proven in Proposition \ref{Thm.Kernels.Notions of equality}.
The implications 2'$\Rrightarrow$1' and 2$\Rrightarrow$1 are obvious.
\\
\\
The implication 1$\Rrightarrow$2. Note that for $(X',0)\norm(X,0)$ the pullback decomposes:
$\nu^*(E)/Torsion=E_1\oplus E_2$. Here $E_i$ is supported on $(X_i,0)$.
As $E$ is $X'/X$ saturated one has: $E=\nu_*(E_1\oplus E_2)\approx E_1\oplus E_2$.

This proves the theorem. Note that in proposition
\ref{Thm.Decomposability.When.Adjoint.Matrix.Belongs.To.Separation.of.Branches} we give also a direct proof
of 1'$\Rrightarrow$2', purely in terms of linear algebra.
\epr
%
%
%
%
%
\section{Some applications}\label{Sec.Applications}
We restrict here to the case of curves, for higher dimensions cf. \S\ref{Sec.Decomposability.Higher.Dim.Case}.
Let $\cM$ be a \mg\ \ dr\ of the plane curve $(C,0)$. Then $\cM$ is decomposable according to the tangential
decomposition (theorem \ref{Thm.Decomposability.For.Tangential.Decomp.Curves}). We study its blocks, each of them is
a \mg\  \dr\ of a curve singularities whose tangent cone has just one line.
\bcor\label{Thm.Classification.of.Max.Gen.dr.for.y.k+x.kl}
Let $(C,0)=\cup(p_\al C_\al)$, where each $C_\al$ is smooth and $T_{(C,0)}=\{x^p_1=0\}$.
Let $\cM$ be a \dr\ of $(X,0)$, \mg\ at the origin.

Then $\cM$ is equivalent to
\beq\bpm f_1&\be_1 x^{n_1}_2&h_{13}(x_2)&..&..&h_{1n}(x_2)
\\0&f_2&\be_2 x^{n_2}_2&h_{24}(x_2)&..&h_{2n}(x_2)
\\..&..&..&..&\\0&..&..&..&0&f_p
\epm\eeq
with  $\be_i\in\{0,1\}$ and  either $h_{ij}(x_2)\equiv0$ or $h_{ij}(x_2)$ a polynomial in $x_2$ such that
$ord_{x_2}(h_{ij})\ge1$ and $deg(h_{ij})<\min(l_i,l_j)$.
\ecor
\bpr 1. The matrix is equivalent to an upper triangular by theorems
\ref{Thm.Decomposability.Max.Gen.Curves.Upper.Block.Triang} and \ref{Thm.Decomposability.Max.Gen.Det.Rep.Multiple.Curve}.
Then by columns subtraction one can kill all the $x_1$ dependent terms in the entries above the diagonal.
\\2. Consider the diagonal $\{(i,i+1)\}$.
Represent each nonzero element $\cM_{i,i+1}(x_2)$ as $x^{n_i}_2\tilde\cM_{i,i+1}$, where
$\tilde\cM_{i,i+1}|_{(0,0)}\neq0$,
i.e. is locally invertible. If $n_i\ge\min(l_i,l_{i+1})$ then by adding the $i'$th column to the column $(i+1)$
and subtracting the row $(i+1)$ from the row $i$ the $x_2$- order can be increased. Continue this process inductively,
 thus killing this entry. Hence, if for some element $\cM_{i,i+1}$ the $x_2$-order is at least $l_i$
 or $l_{i+1}$  the element can be just set to zero.  The remaining non-zero elements $x^{n_i}\tilde\cM_{i,i+1}$ are set to $x^{n_i}$
 by the conjugation $\cM\to U^{-1}\cM U$ with
\beq
U=\bpm \prodl_{i\ge1}\tilde\cM_{i,i+1}&0&..&0\\0&\prodl_{i\ge2}\tilde\cM_{i,i+1}&..&0\\0&..&..&0&\tilde\cM_{k-1,k}\epm
\eeq
Regarding the remaining entries $h_{ij}(x)$ with $j-i\ge2$, bring them to the needed form diagonal-by-diagonal.
This is done again by the standard procedure: add $y+x^{l_i}$, subtract $y+x^{l_j}$ etc.
\epr

\bex
$\bullet$ Any \mg\ \dr\ of $y(y+x^{l_1})(y-x^{l_2})$ is equivalent to either:
\beq
\bpm y+x^{l_1}&x^{n_1}&h(x)\\0&y&x^{n_2}\\0&0&y-x^{l_2}\epm,\quad n_i< l_i,\quad 1\le ord(h(x))<\min(n_1,n_2)
\text{ or }h(x)\equiv0
\eeq
or
\beq
\bpm y+x^{l_1}&0&x^{n_1}\\0&y&x^{n_2}\\0&0&y-x^{l_2}\epm,\quad
\bpm y+x^{l_1}&x^{n_1}&x^{n_2}\\0&y&0\\0&0&y-x^{l_2}\epm,\quad
\bpm y+x^{l_1}&0&x^{n_1}\\0&y&0\\0&0&y-x^{l_2}\epm
\eeq
\li For $(C,0)=\{y(y^2-x^{2l+1})=0\}$ any \mg\ \dr\ is equivalent to
\beq
\bpm y&p_1(x)&p_2(x)\\0&y&x^{2l+1-m}\\0&x^m&y\epm,\quad p_1(0)=0=p_2(0),\quad \deg(p_1(x))<m,\quad \deg(p_2(x))<2l+1-m
\eeq
Here the pair of polynomials $(p_1(x),p_2(x))$ is determined up to scaling. Hence the space of \mg\ representations
is parameterized by $H^0(\cO_{\P^1}(m-2))\times H^0(\cO_{\P^1}(2l-m))/\sim$, where $H^0(\cO_{\P^1}(j))$ is the space
of homogeneous polynomials in two variables of degree $j$ and the equivalence relation is the scaling.
\eex

\end{document}